\documentclass[a4paper,12pt]{article}

\usepackage{amsmath,amssymb,latexsym,amsfonts,amsthm}
\usepackage{mathrsfs}

\theoremstyle{plain}

\theoremstyle{definition}

\newcommand{\bD}{\ensuremath{\mathbb{D}}}
\newcommand{\bE}{\ensuremath{\mathbb{E}}}

\newcommand{\bR}{\ensuremath{\mathbb{R}}}

\newcommand{\cB}{\ensuremath{\mathcal{B}}}

\newcommand{\cE}{\ensuremath{\mathcal{E}}}
\newcommand{\cF}{\ensuremath{\mathcal{F}}}

\newcommand{\vm}{\ensuremath{\mbox{{\boldmath $m$}}}}
\newcommand{\vn}{\ensuremath{\mbox{{\boldmath $n$}}}}

\renewcommand{\Re}{{\rm Re} \ }

\newcommand{\eps}{\ensuremath{\varepsilon}}
\newcommand{\e}{{\rm e}}
\renewcommand{\d}{{\rm d}}

\newcommand{\law}{\stackrel{{\rm law}}{=}}

\renewcommand{\hat}{\widehat}

\newcommand{\Supp}{\mathop{\rm Supp}}

\newcommand{\abra}[1]{\left| #1 \right|}
\newcommand{\cbra}[1]{\left( #1 \right)}
\newcommand{\kbra}[1]{\left\{ #1 \right\}}
\newcommand{\ebra}[1]{\left[ #1 \right]}

\newcommand{\n}{\nonumber}

\numberwithin{equation}{section}

\newcounter{No}

\newcounter{Ci}[subsection]

\makeatletter
\renewcommand\section{\@startsection {section}{1}{\z@}%
                                   {-3.5ex \@plus -1ex \@minus -.2ex}%
                                   {2.3ex \@plus.2ex}%
                                   {\normalfont\large\bf}}
\makeatother

\newtheorem{theorem}{Theorem}[section]
\newtheorem{proposition}{Proposition}[section]
\newtheorem{corollary}{Corollary}[section]
\newtheorem{lemma}{Lemma}[section]
\newtheorem{example}{Example}[section]
\newtheorem*{remark}{Remark}

\newcommand{\harm}{h}
\renewcommand{\harm}{h}
\newcommand{\sym}{h,s}

\begin{document}

\title{Excursions away from a regular point 
for one-dimensional symmetric L\'evy processes without Gaussian part\thanks{
The research of the author is supported by KAKENHI (20740060)}
}
%\subtitle{Do you have a subtitle?\\ If so, write it here}

%\titlerunning{Excursions for L\'evy processes}        % if too long for running head

\author{Kouji Yano\footnote{Graduate School of Science, Kobe University} 
}

\maketitle

\begin{abstract}
The characteristic measure of excursions away from a regular point is studied 
for a class of symmetric L\'evy processes without Gaussian part. 
It is proved that 
the harmonic transform of the killed process enjoys Feller property. 
The result is applied to prove extremeness of the excursion measure 
and to prove several sample path behaviors of the excursion and the $ h $-path processes. 

\


\end{abstract}

\section{Introduction}

It\^o \cite{MR0402949} has proved that the point process of excursions 
away from a regular point for a strong Markov process is Poisson. 
Its characteristic measure 
will be simply called the {\em excursion measure}. 
It\^o's theorem shows that, 
for a given minimal process, 
there is a one-to-one correspondence 
between an excursion measure and a strong Markov extension of the minimal process. 
In the same paper, he established the integral representation formula 
on the convex set of the normalized excursion measures of strong Markov extensions 
of the minimal process. 

In the present paper, we study several properties 
of the excursion measure and the $ h $-path process 
for one-dimensional symmetric L\'evy processes. 
Under certain assumptions which imply no Gaussian part, 
we prove extremeness of the excursion measure 
and several sample path behaviors of the excursion and the $ h $-path processes. 
Our study is motivated by a recent study of Yano--Yano--Yor \cite{YYY} 
about penalization problems for one-dimensional 
symmetric $ \alpha $-stable processes of index $ 1<\alpha \le 2 $.

Let $ \{ (X_t:t \ge 0),(\cF_t:t \ge 0),(P_x:x \in \bR) \} $ 
denote the canonical representation of the one-dimensional Brownian motion. 
Let $ \vn^{\rm B} $ stand for the {\em Brownian excursion measure}. 
Then the measure $ \vn^{\rm B} $ is represented as 
\begin{align}
\vn^{\rm B} = \frac{1}{2} \vn^{+} + \frac{1}{2} \vn^{-} 
\label{eq: exc meas of skew BM}
\end{align}
where $ \vn^{+} $ (resp. $ \vn^{-} $) stands for the excursion measure 
for (resp. the negative of) the reflecting Brownian motion. 
We remark that the formula \eqref{eq: exc meas of skew BM} is a special case 
of {\em It\^o's integral representation formula} 
(\cite[Theorem 7.1]{MR0402949}; see also \cite[Section V.6]{MR1138461}): 
Let $ \cE $ denote the convex cone 
of the excursion measures for non-trivial strong Markov extensions 
of the killed process. 
Let $ \cE_1 $ denote the convex subset of $ \cE $ 
whose elements $ \mu $ are {\em normalized} in the sense that $ \mu[1-\e^{-\zeta}] = 1 $ 
where $ \zeta $ stands for the lifetime. 
Then any given element $ \mu_0 $ of the set $ \cE_1 $ is represented as 
\begin{align}
\mu_0(\cdot) = \int_{\displaystyle {\rm ex} (\cE_1)} \mu(\cdot) \Pi(\d \mu) 
\n
\end{align}
for some probability measure $ \Pi $ 
on the set $ {\rm ex} (\cE_1) $ of extreme points of $ \cE_1 $. 
We say that $ \vn \in \cE $ is an {\em extreme direction} 
if $ \vn $ is proportional to an extreme point of $ \cE_1 $, 
i.e., 
\begin{align}
\frac{\vn(\cdot)}{\vn[1-\e^{-\zeta}]} \in {\rm ex} (\cE_1). 
\n
\end{align}

Let $ \{ P^{+}_x:x \ge 0 \} $ (resp. $ \{ P^{-}_x:x \le 0 \} $) 
denote the law of (resp. the negative of) the three-dimensional Bessel process. 
Let $ P^{+,-}_x $ denote the law of the symmetrized three-dimensional Bessel process, 
i.e., 
$ P^{+,-}_x = P^{+}_x $ for $ x>0 $, 
$ P^{+,-}_x = P^{-}_x $ for $ x<0 $ 
and $ P^{+,-}_0 = \frac{1}{2} P^{+}_0 + \frac{1}{2} P^{-}_0 $. 
The process $ \{ X_{\cdot},P^{+,-}_{\cdot} \} $ is 
the {\em harmonic transform} or {\em $ h $-path process} 
with respect to the harmonic function $ h(x)=|x| $ 
of the killed Brownian motion $ \{ P^0_x:x \in \bR \setminus \{ 0 \} \} $ 
in the sense that $ \d P^{+,-}_x|_{\cF_t} = \frac{|X_t|}{|x|} \d P^0_x|_{\cF_t} $ 
for any $ x \in \bR \setminus \{ 0 \} $. 
The law $ P^{+,-}_0 $ 
is related to the excursion measure $ \vn^{\rm B} $ 
in the following {\em Imhof relation} 
(see, e.g., \cite[Exercise XII.4.18]{MR1725357}): 
\begin{align}
\d P^{+,-}_0|_{\cF_t} = \frac{|X_t|}{\vn^{\rm B}[|X_t|]} \d \vn^{\rm B}|_{\cF_t} 
, \qquad t>0 . 
\n
\end{align}
Here we remark that the law 
$ \frac{1_{\{ \zeta>t \}}}{\vn^{\rm B}(\zeta>t)} \d \vn^{\rm B}|_{\cF_t} $ 
is nothing but the law of the {\em Brownian meander}. 
Moreover, the excursion measure $ \vn^{\rm B} $ admits 
the following {\em lifetime disintegration formula} 
(see \cite[Section III.4.3]{MR1011252} and \cite[Theorem XII.4.2]{MR1725357}): 
\begin{align}
\vn^{\rm B}(\cdot) = \int_0^{\infty } P^{+,-}_0(\cdot|X_{t-}=0) \frac{\d t}{\sqrt{2 \pi t^3}} . 
\label{eq: LDF of BM}
\end{align}
In other words, under $ \vn^{\rm B} $, the law of the lifetime $ \zeta $ 
is $ \d t/\sqrt{2 \pi t^3} $ and 
the conditional law of the excursion process given $ \zeta=t $ 
is $ P^{\harm}_0(\cdot|X_{t-}=0) $. 

We point out the following three facts: 
\\ \quad (i) 
The excursion measure $ \vn^{\rm B} $, 
being represented as \eqref{eq: exc meas of skew BM}, 
is {\em not} an extreme direction; 
\\ \quad (ii) 
The {\em germ $ \sigma $-field} $ \cF_{0+} = \cap_{\eps>0} \cF_{\eps} $ 
is not trivial under $ \vn^{\rm B} $; 
in fact, the set 
$ A = \{ \exists t>0 \ \text{such that} \ \forall s \le t, \ X_s \ge 0 \} $ 
belongs to $ \cF_{0+} $, 
but {\em neither $ A $ nor} $ A^c $ is $ \vn^{\rm B} $-null; 
\\ \quad (iii) 
The semigroup $ \{ T^{+,-}_t:t \ge 0 \} $ corresponding to the $ h $-path process 
does {\em not} enjoy Feller property; 
in fact, for a positive continuous function $ f $ with compact support in $ [0,\infty ) $, we have 
\begin{align}
\lim_{x \to 0+} T^{+,-}_tf(x) = P^{+}_0[f(X_t)] >0 
\qquad \text{while} \qquad 
\lim_{x \to 0-} T^{+,-}_tf(x) = 0 . 
\n
\end{align}

\subsection*{Main theorems}

Let us state the main theorems of the present paper. 
Let $ \{ P_x:x \in \bR \} $ denote the law of a one-dimensional L\'evy process. 
We assume that the following conditions are satisfied: 
\begin{quote}
{\bf (A0)} 
The process is symmetric; 
\\
{\bf (A1)} 
The origin is regular for itself; 
\\
{\bf (A2)} 
The process is not a compound Poisson. 
\end{quote}
The set of these conditions {\bf (A0)}-{\bf (A2)} will be denoted simply by {\bf (A)}. 
Then it holds (see Section \ref{sec: KB}) that 
the {\em L\'evy--Khintchine exponent} 
\begin{align}
\theta(\lambda) = v \lambda^2 + 2 \int_{(0,\infty )} (1-\cos \lambda x) \nu(\d x) 
, \qquad \lambda \in \bR 
\label{def theta}
\end{align}
with the {\em Gaussian coefficient} $ v $ and the {\em L\'evy measure} $ \nu(\d x) $ satisfies 
\begin{align}
\int_0^{\infty } \frac{\min \{ \lambda^2 , 1 \}}{\theta(\lambda)} \d \lambda < \infty . 
\label{eq: theta lambda ibility}
\end{align}
Hence there exists a continuous density $ u_q(x) $ of the resolvent kernel given by 
\begin{align}
u_q(x) = \frac{1}{\pi} \int_0^{\infty } \frac{\cos \lambda x}{q + \theta(\lambda)} \d \lambda 
\n
\end{align}
and, in addition, the following function is well-defined (see also \cite[Lemma 1]{SY}): 
\begin{align}
h(x) 
= \lim_{q \to 0+} \kbra{ u_q(0)-u_q(x) } 
= \frac{1}{\pi} \int_0^{\infty } \frac{1-\cos \lambda x}{\theta(\lambda)} \d \lambda . 
\label{hx}
\end{align}
Remark that the process is recurrent or transient according as 
\begin{align}
\kappa = \lim_{q \to 0+} \frac{1}{u_q(0)} 
= \kbra{ \frac{1}{\pi} \int_0^{\infty } \frac{1}{\theta(\lambda)} \d \lambda }^{-1} 
\label{eq: kappa}
\end{align}
is zero or positive 
(see Section \ref{sec: rec and trans}). 

Let $ T_{\{ x \}} $ denote the first hitting time of $ x \in \bR $. 
Let $ L(t,x) $ denote the local time process. 
We denote by $ \{ P^0_x:x \in \bR \} $ the law of the process 
killed upon hitting the origin, 
which we simply call the {\em killed process} in short. 
We denote by $ \vn $ the excursion measure. 

\begin{theorem} \label{thm: harmonic}
Suppose that the condition {\bf (A)} is satisfied. Then 
the function $ h(x) $ is invariant excessive with respect to the killed process, i.e., 
\begin{align}
P^0_x[h(X_t)] = h(x) 
, \qquad t>0 , \ x \in \bR \setminus \{ 0 \} . 
\n
\end{align}
\end{theorem}

\begin{remark}
The function $ h(x) $ above is {\em harmonic} 
in the sense of \cite[Definition 4.3.2]{MR2152573} 
where the Laplacian is replaced by the generator. 
\end{remark}

\begin{theorem} \label{thm: vn hXt = 1}
Suppose that the condition {\bf (A)} is satisfied. Then 
\begin{align}
\vn [h(X_t)] = 1 
, \qquad t > 0 . 
\n
\end{align}
\end{theorem}

Theorems \ref{thm: harmonic} and \ref{thm: vn hXt = 1} 
will be proved in Section \ref{sec: h-func}. 

We introduce the $ h $-path process $ \{ P^{\harm}_x:x \in \bR \} $ 
as the law on the canonical space such that, for any $ t>0 $, 
\begin{align}
\d P^{\harm}_x |_{\cF_t} =& 
\frac{h(X_t)}{h(x)} \d P^0_x |_{\cF_t} 
, \qquad x \in \bR \setminus \{ 0 \} , 
\label{eq: RN density1}
\\
\d P^{\harm}_0 |_{\cF_t} =& 
h(X_t) \d \vn |_{\cF_t} 
, \qquad \ \ x=0 . 
\label{eq: RN density2}
\end{align}
Note that such a family of probability laws on $ \bD $ exists uniquely, 
because the identities \eqref{eq: RN density1} and \eqref{eq: RN density2} 
induce a consistent family of probability laws on $ \bD $ 
by the Markov properties 
of the killed process and the excursion process. 
Let us denote the corresponding semigroup by $ \{ T^{\harm}_t:t \ge 0 \} $. 

\begin{theorem}[Lifetime disintegration formula] \label{thm: disinteg by lifetime}
Suppose that the condition {\bf (A)} is satisfied. Then 
there exists a positive completely monotone function $ \rho(t) $ such that 
\begin{align}
\vn(\cdot) 
= \int_0^{\infty } P^{\harm}_0(\cdot | X_{t-}=0) \rho(t) \d t 
+ \kappa P^{\harm}_0(\cdot) . 
\label{eq: lifetime disinteg}
\end{align}
In other words, 
\\
\quad {\rm (i-a)} 
$ \vn(\zeta \in \d t) = \rho(t) \d t $ on $ (0,\infty ) $; 
\\
\quad {\rm (i-b)} 
$ \vn(\cdot | \zeta =t) = P^{\harm}_0(\cdot | X_{t-}=0) $ for $ 0<t<\infty $; 
\\
\quad {\rm (ii-a)} 
$ \vn(\zeta = \infty ) = \kappa $; 
\\
\quad {\rm (ii-b)} 
$ \vn(\cdot | \zeta =\infty ) = P^{\harm}_0(\cdot) $. 
\end{theorem}

Theorem \ref{thm: disinteg by lifetime} 
will be proved 
in Section \ref{sec: proof of Lifetime disinteg}. 
We must remark that Getoor--Sharpe \cite[Theorem 7.6]{MR673804} 
have proved that the excursion measure for quite general Markov processes 
admits lifetime disintegration formula 
with a certain bridge process as its conditional distribution. 
(Note that $ P^{0,t,0} $ in \cite{MR673804} corresponds to our $ P^h_0(\cdot | X_{t-}=0) $, 
$ \eta(t,0,0) $ to $ \rho(t) $, 
$ q^*(t,0,x) $ to $ \rho(t,x) $, 
and $ q(t,x,y) $ to $ p^0_t(x,y) $.) 
Theorem \ref{thm: disinteg by lifetime} asserts that, 
in this particular case, 
the conditional distribution is given by 
the bridge process of the $ h $-path process. 
In the same way as Theorem \ref{thm: disinteg by lifetime}, we may prove that, 
for any $ x \in \bR \setminus \{ 0 \} $, 
\begin{align}
P^0_x(\cdot) = \int_0^{\infty } P^h_x(\cdot|X_{t-}=0) P^0_x(\zeta \in \d t) 
+ \kappa h(x) P^h_x(\cdot) . 
\label{eq: LD for P0x}
\end{align}
In particular, we have $ P^0_x(\cdot | \zeta=t) = P^h_x(\cdot|X_{t-}=0) $. 

Although it seems superflous, 
we need the following extra assumption for some technical reason: 
\begin{quote}
{\bf (T)} 
The function $ \theta(\lambda) $ is non-decreasing in $ \lambda>\lambda_0 $ 
for some $ \lambda_0>0 $. 
\end{quote}
The following theorem asserts that the $ h $-path process is transient. 

\begin{theorem} \label{thm: trans of h-path}
Suppose that the conditions {\bf (A)} and {\bf (T)} are satisfied. Then 
\begin{align}
P^{\harm}_x \cbra{ \lim_{t \to \infty } |X_t| = \infty } = 1 
, \qquad x \in \bR . 
\n
\end{align}
\end{theorem}

Theorem \ref{thm: trans of h-path} will be proved in Section \ref{sec: transience}.

We need the following assumption: 
\begin{quote}
{\bf (B)} 
$ \displaystyle \lim_{x \to 0} \frac{x}{h(x)} = 0 $. 
\end{quote}
From the assumption {\bf (B)} 
it follows that the Gaussian coefficient $ v $ is zero (see Lemma \ref{lem: excludes}). 
Now let us state our main theorem. 

\begin{theorem} \label{thm: Feller}
Suppose that the conditions {\bf (A)}, {\bf (B)} and {\bf (T)} are satisfied. Then 
the semigroup $ \{ T^{\harm}_t:t \ge 0 \} $ enjoys Feller property. 
\end{theorem}

As applications of Theorem \ref{thm: Feller}, we obtain 

\begin{corollary}[Extremeness property] \label{cor: extreme}
Suppose that the conditions {\bf (A)}, {\bf (B)} and {\bf (T)} are satisfied. Then 
the excursion measure $ \vn $ is an extreme direction. 
\end{corollary}

\begin{corollary}[Oscillatory entrance property] \label{cor: change signs of h-path}
Suppose that the conditions {\bf (A)}, {\bf (B)} and {\bf (T)} are satisfied. Then 
the excursion process enters oscillatingly, 
i.e., 
\begin{align}
\vn \cbra{ \kbra{ \text{$ \exists \{ t_n \} $ with $ t_n \searrow 0 $ such that 
$ \forall n $, $ X_{t_n} X_{t_{n+1}} < 0 $} }^c } = 0 . 
\n
\end{align}
\end{corollary}

The following results are concerned about sample path behaviors 
of the $ h $-path processes. 

\begin{corollary}[Oscillatory entrance property] \label{cor: change signs of h-path-2}
Suppose that the conditions {\bf (A)}, {\bf (B)} and {\bf (T)} are satisfied. Then 
the $ h $-path process enters oscillatingly, 
i.e., 
\begin{align}
P^{\harm}_0 \cbra{ \text{$ \exists \{ t_n \} $ with $ t_n \searrow 0 $ such that 
$ \forall n $, $ X_{t_n} X_{t_{n+1}} < 0 $} } = 1 . 
\n
\end{align}
\end{corollary}

\begin{corollary}[Oscillatory property in the long time] \label{cor: oscillation of stable}
Suppose that 
the process $ \{ (X_t),(P_x) \} $ is a symmetric stable process of index $ 1<\alpha <2 $. 
Then 
\begin{align}
P^{\harm}_0 \cbra{ \limsup_{t \to \infty } X_t 
= \limsup_{t \to \infty } (-X_t) = \infty } = 1 . 
\n
\end{align}
\end{corollary}

Theorem \ref{thm: Feller} 
and Corollaries \ref{cor: extreme}, 
\ref{cor: change signs of h-path}, 
\ref{cor: change signs of h-path-2} 
and \ref{cor: oscillation of stable} will be proved 
in Section \ref{sec: Feller}. 
Here we briefly sketch how to prove Corollary \ref{cor: extreme}, 
provided that Theorem \ref{thm: Feller} is proved, 
as follows:  
\begin{align}
& \text{ The semigroup $ \{ T^{\harm}_t: t \ge 0 \} $ enjoys Feller property } 
\\
\stackrel{\text{Prop. \ref{prop: germ}}}{\Longrightarrow} & 
\text{ The germ $ \sigma $-field $ \cF_{0+} $ is trivial under $ P^{\harm}_0 $ } 
\\
\stackrel{\text{eq. \eqref{eq: RN density2} }}{\Longleftrightarrow} & 
\text{ The germ $ \sigma $-field $ \cF_{0+} $ is trivial under $ \vn $ } 
\\
\stackrel{\text{Thm. \ref{thm: extremeness}}}{\Longleftrightarrow} & 
\text{ The excursion measure $ \vn $ is an extreme direction. } 
\n
\end{align}

\begin{example}
When it is a stable process, 
the process satisfies the condition {\bf (A)} 
if and only if it is a symmetric $ \alpha $-stable process of index $ 1<\alpha \le 2 $. 
Up to multiplicative constant, 
we have 
\begin{align}
\theta(\lambda)=|\lambda|^{\alpha } 
\qquad \text{for some $ 1 < \alpha \le 2 $}. 
\n
\end{align}
Hence the condition {\bf (T)} is automatically satisfied. 
The harmonic function is given by 
\begin{align}
h(x) = C(\alpha ) |x|^{\alpha -1} 
, \qquad x \in \bR 
\n
\end{align}
where 
\begin{align}
C(\alpha ) = \frac{1}{\pi} \int_0^{\infty } \frac{1-\cos \lambda}{\lambda^{\alpha }} \d \lambda 
= 
\begin{cases}
\frac{\Gamma (2-\alpha )}{\pi (\alpha -1)} \sin \frac{\alpha \pi}{2} 
\quad & \text{for} \ 1<\alpha <2, \\
\frac{1}{2} 
\quad & \text{for} \ \alpha =2. 
\end{cases}
\n
\end{align}
The density $ \rho(t) $ and the constant $ \kappa $ 
in the formula \eqref{eq: lifetime disinteg} are given by 
\begin{align}
\rho(t) = \frac{(\alpha -1) \pi}{\Gamma(1-1/\alpha ) \Gamma(1/\alpha )^2} 
t^{\frac{1}{\alpha } -2 } 
\qquad \text{and} \qquad 
\kappa = 0 . 
\n
\end{align}
Note that the condition {\bf (B)} is satisfied if and only if $ 1<\alpha <2 $. 
\end{example}

The organization of the present paper is as follows. 
In Section \ref{sec: prel}, we recall several preliminary facts 
about one-dimensional L\'evy processes in general settings. 
In Section \ref{sec: prel sym}, we recall several preliminary facts 
assuming that the process is symmetric. 
In Section \ref{sec: h-func}, 
we prove harmonicity of $ h(x) $ and study its properties. 
In Section \ref{sec: h-trans}, 
we prove the lifetime disintegration formula for the excursion measure. 
In Section \ref{sec: key lemmas}, 
we prove several lemmas for later use. 
The transience of the $ h $-path process will be proved there. 
Section \ref{sec: Feller} 
is devoted to the proof of Feller property of the semigroup 
corresponding to the $ h $-path process. 
The extremeness property of the excursion measure 
and the sample path behaviors of the excursion and the $ h $-path processes 
will be proved in this section.

\subsection*{Remarks}

\begin{remark}
Contrary to that for diffusion processes, 
boundary problem for Markov processes with jumps is extremely difficult 
because of non-locality. 
We must remark that 
Chen--Fukushima--Ying \cite[Sections 4 and 5]{CFY} (see also Fukushima--Tanaka \cite{MR2139028}) 
have proved under quite a general assumption that 
there exists a unique extension of the minimal process which conserves a given weak duality. 
Thanks to this striking result, we know that, at least in the settings of the present paper, 
the symmetric extension of the minimal process is unique. 
\end{remark}

\begin{remark}
Based on a kind of Martin boundary argument (see Lemma \ref{lem: u dagger is conti at 0}), 
Ikeda--Watanabe \cite{MR0451425} (see also Takada \cite{MR0370784}) 
have studied sample path behaviors before hitting the origin 
for (possibly non-symmetric) one-dimensional L\'evy processes. 
With the help of their results, we can give another proof 
of a special case of Corollary \ref{cor: change signs of h-path} 
without using Feller property of the $ h $-path process (see Section \ref{sec: Ikeda-Watanabe}). 
\end{remark}

\begin{remark}
Several aspects of excursion measures have been extensively studied 
for Brownian motions 
(see, e.g., \cite[Section III.4.3]{MR1011252} and \cite[Chap. XII]{MR1725357}), 
for diffusion processes 
(see, e.g., \cite{MR2295612}; see also \cite{MR2384481}), 
for reflected L\'evy processes 
(see, e.g., \cite[Chap.VI]{MR1406564}) 
and for spectrally one-sided L\'evy processes 
(see, e.g., \cite[Chap.VII]{MR1406564}). 
In particular, 
the lifetime disintegration formula 
for one-dimensional diffusion processes 
can be found in Pitman--Yor \cite{MR656509,MR1439532} and Yano \cite{MR2266999}. 

For symmetric L\'evy processes, however, we cannot find any literature 
about the excursion measure, 
except general theories 
and Fitzsimmons--Getoor \cite{MR1341555} (see also Yano--Yano \cite{YY}) 
who have studied the law of the time spent on the positive side 
by the conditional process $ \{ (X_t:0 \le t \le T),\vn(\cdot|\zeta=T) \} $. 
Note that Theorem \ref{thm: disinteg by lifetime} asserts that 
the conditional law is given by the bridge $ P^{\harm}_0(\cdot|X_{T-}=0) $ 
of the $ h $-path process. 
\end{remark}

\begin{remark}
In order to study regularity of the local time $ L(t,x) $, 
Barlow \cite{MR775850} (see also \cite[Section V.3]{MR1406564}) has introduced 
the following function: 
\begin{align}
h_B(x) 
= P_x[L(T_{\{ 0 \}},x)] 
= \lim_{q \to 0+} \kbra{ u_q(0) - \frac{u_q(x) u_q(-x)}{u_q(0)} } . 
\n
\end{align}
The function $ h_B $ is related to our $ h $ as follows: 
\begin{align}
h_B(x) 
= 2h(x) - \kappa h(x)^2 . 
\n
\end{align}
In particular, $ h_B = 2h $ in the recurrent case. 
\end{remark}

\begin{remark}
Salminen--Yor \cite{SY} have obtained the following Tanaka formula: 
\begin{align}
h(X_t-x) = h(x) + N_t^x + c_1 L(t,x) 
\n
\end{align}
where $ N_t^x $ is a $ P_0 $-martingale 
and $ c_1 $ is some constant. 
\end{remark}

\begin{remark}
The $ h $-path process $ (X_{\cdot},P^{\harm}_{\cdot}) $ 
is considered to be the process $ (X_{\cdot},P_{\cdot}) $ 
conditioned never to hit the origin. 
In fact, in the case of symmetric $ \alpha $-stable process of index $ 1<\alpha \le 2 $, 
Yano--Yano--Yor \cite{YYY} proved that 
\begin{align}
\lim_{t \to \infty } P_x[Z_s|T_{\{ 0 \}}>t] = P^{\harm}_x[Z_s] 
\n
\end{align}
for all non-negative $ \cF_s $-measurable functional $ Z_s $. 
\end{remark}

\section{Preliminary facts: general case} \label{sec: prel}

Let $ \bD $ denote the set of c\`adl\`ag paths 
$ w:[0,\infty ) \to \bR \cup \{ \Delta \} $ such that 
$ w(t)=\Delta $ for all $ t \ge \zeta(w) $ where 
\begin{align}
\zeta(w) = \inf \{ t \ge 0 : w(t)=\Delta \} . 
\n
\end{align}
Here the topology of $ \bR \cup \{ \Delta \} $ is the one-point compactification of $ \bR $. 
The point $ \Delta $ is called the {\em cemetery} 
and $ \zeta(w) $ is called the {\em lifetime} of a path $ w \in \bD $. 
The space $ \bD $ is equipped with Skorokhod topology. 
Let $ \cB_{+,b}(\bR) $ denote the set of measurable functions 
which are non-negative or bounded. 
For $ f \in \cB_{+,b}(\bR) $, we define $ \| f \| = \sup_{x \in \bR} |f(x)| $. 
Let $ C_0(\bR) $ denote the set of continuous functions 
which vanish at infinity, i.e., $ \lim_{|x| \to \infty } f(x) = 0 $. 

Let $ (X_t:t \ge 0) $ denote the coordinate process: $ X_t(w)=w(t) $, $ t \ge 0 $. 
Let $ (\cF_t: t \ge 0) $ denote the natural filtration: 
$ \cF_t = \sigma(X_s: s \le t) $. 
Let $ (P_x: x \in \bR) $ be the law 
of a one-dimensional L\'evy process 
on the canonical space $ \bD $. 
Throughout this section, 
we do {\em not} suppose that $ \{ (X_t),(P_x) \} $ is symmetric. 
The corresponding semigroup and resolvent operator will be denoted by 
\begin{align}
T_tf(x) = P_x[f(X_t)] 
\qquad \text{and} \qquad 
U_qf(x) = \int_0^{\infty } \e^{-qt} T_tf(x) \d t , 
\label{2.2}
\end{align}
respectively. 
It is well-known that 
\begin{align}
P_0[\e^{i \lambda X_t}] = \e^{t \psi(\lambda)} 
, \qquad \lambda \in \bR 
\n
\end{align}
where the L\'evy--Khintchine exponent $ \psi(\lambda) $ is given by 
\begin{align}
\psi(\lambda) = - v \lambda^2 + i a \lambda 
+ \int_{\bR} \cbra{ \e^{i \lambda x} - 1 - \frac{i \lambda x}{1+x^2} } \nu(\d x) 
, \qquad \lambda \in \bR 
\n
\end{align}
for some $ v \ge 0 $, $ a \in \bR $ 
and some positive Radon measure $ \nu $ on $ \bR $ such that 
\begin{align}
\int_{\bR} \min \{ x^2,1 \} \nu(\d x) < \infty . 
\n
\end{align}
Set 
\begin{align}
\theta(\lambda) 
= - \Re \psi(\lambda) 
= v \lambda^2 + \int_{\bR} \cbra{ 1-\cos \lambda x } \nu(\d x) 
, \qquad \lambda \in \bR . 
\n
\end{align}

\subsection{Germ triviality}

It is obvious that 
the semigroup $ \{ T_t:t \ge 0 \} $ enjoys the Feller property; 
in particular, 
\begin{align}
T_t C_0(\bR) \subset C_0(\bR) 
, \qquad t \ge 0 . 
\label{eq: Tt C0 subset C0}
\end{align}

\begin{proposition}[Blumenthal \cite{MR0088102}] \label{prop: germ}
For any $ x \in \bR $, 
the germ $ \sigma $-field $ \cF_{0+} $ is trivial under $ P_x $. 
\end{proposition}

The following proof says that 
{\em Feller property implies germ triviality}. 

\begin{proof}
Let $ A \in \cF_{0+} $, $ f \in C_0(\bR) $ 
and $ t,\eps>0 $. 
By the Markov property, we have 
\begin{align}
P_x[1_A f(X_{t+\eps})] 
= P_x[1_A T_tf(X_{\eps})] . 
\n
\end{align}
Now let $ \eps $ tend to $ 0+ $. 
By the continuity of $ T_tf $, 
by right-continuity of paths 
and by the dominated convergence theorem, 
we obtain 
\begin{align}
P_x[1_A f(X_t)] 
= P_x(A) T_tf(x) 
= P_x(A) P_x[f(X_t)] . 
\n
\end{align}
The rest of the proof is a standard argument. 
\qed
\end{proof}

\subsection{The condition for the origin to be regular for itself} \label{sec: cond A}

Recall the following conditions: 
\begin{quote}
{\bf (A1)} 
The origin is regular for itself; 
\\
{\bf (A2)} 
The process is not a compound Poisson. 
\end{quote}

\begin{theorem}[Kesten \cite{MR0272059} and Bretagnolle \cite{MR0368175}] \label{thm: Kesten}
The conditions {\bf (A1)}-{\bf (A2)} are satisfied if and only if 
\begin{align}
\int_{\bR} \Re \frac{1}{q-\psi(\lambda)} \d \lambda < \infty 
, \qquad q>0 
\label{eq: KB}
\end{align}
and 
\begin{align}
\text{either} \qquad 
v>0 
\qquad \text{or} \qquad 
\int_{(-1,1)} |x| \nu(\d x) = \infty . 
\label{eq: KB2}
\end{align}
\end{theorem}

In what follows, we suppose that the conditions {\bf (A1)}-{\bf (A2)} are satisfied. 
By Fourier inversion, we see that the function 
\begin{align}
p_t(x) = \frac{1}{2\pi} \int_{\bR} \Re \e^{- i \lambda x + t \psi(\lambda)} \d \lambda 
, \qquad t>0 , \ x \in \bR 
\label{ptx}
\end{align}
is a continuous density of the transition probability: 
\begin{align}
P_x(X_t \in A) = \int_A p_t(y-x) \d y 
\n
\end{align}
for $ t>0 $, $ x \in \bR $ and $ A \in \cB(\bR) $. 
The Laplace transform of $ p_t(x) $: 
\begin{align}
u_q(x) 
= \int_0^{\infty } \e^{-qt} p_t(x) \d t 
= \frac{1}{2 \pi} \int_{\bR} \Re \frac{\e^{- i \lambda x}}{q - \psi(\lambda)} \d \lambda 
, \qquad q>0 , \ x \in \bR 
\label{uqx}
\end{align}
is a continuous density of the resolvent kernel: 
\begin{align}
P_x \ebra{ \int_0^{\infty } \e^{-qt} 1_A(X_t) \d t } 
= \int_A u_q(y-x) \d y 
\n
\end{align}
for $ q>0 $, $ x \in \bR $ and  $ A \in \cB(\bR) $. 
We note that the resolvent equation 
\begin{align}
U_q - U_r + (q-r) U_q U_r = 0 
\n
\end{align}
implies that 
\begin{align}
\int u_q(y-x) u_r(z-y) \d y = \frac{1}{q-r} \kbra{ u_r(z-x) - u_q(z-x) } 
\label{eq: resolv eq}
\end{align}
for all $ q,r>0 $ with $ q \neq r $ and $ z,x \in \bR $.

\begin{theorem} \label{thm: hitting}
Suppose that the conditions {\bf (A1)}-{\bf (A2)} are satisfied. Then 
\begin{align}
P_x \ebra{ \e^{-q T_{\{ 0 \}}} } = \frac{u_q(-x)}{u_q(0)} 
, \qquad 
x \in \bR, \ q>0 . 
\label{hitting Lap trans}
\end{align}
\end{theorem}

Theorem \ref{thm: hitting} can be proved 
via {\em Hunt's switching identity}, 
which is based on the duality between $ (X_t) $ and $ (-X_t) $. 
The proof can be found, e.g., in \cite[pp. 64]{MR1406564}, 
and so we omit it.

\subsection{Killed process}

We define the killed process as 
\begin{align}
X^0_t = 
\begin{cases}
X_t    \ & t < T_{\{ 0 \}} , \\
\Delta \ & t \ge T_{\{ 0 \}} . 
\end{cases}
\label{eq: killed process}
\end{align}
The laws on the space $ \bD $ 
of the killed process $ (X^0_t:t \ge 0) $ under $ \{ P_x:x \in \bR \setminus \{ 0 \} \} $ 
will be denoted by 
$ (P^0_x:x \in \bR \setminus \{ 0 \}) $. 
The corresponding semigroup and resolvent operator will be denoted by 
$ T^0_t $ and $ U^0_q $, respectively, in the same way as \eqref{2.2}. 

Set 
\begin{align}
p^0_t(x,y) = p_t(y-x) - \int_0^t p_{t-s}(y) P_x(T_{\{ 0 \}} \in \d s) 
\label{p0txy}
\end{align}
and 
\begin{align}
u^0_q(x,y) 
= \int_0^{\infty } \e^{-qt} p^0_t(x,y) \d t 
= u_q(y-x) - \frac{u_q(-x) u_q(y)}{u_q(0)} 
\label{uq0xy}
\end{align}
for $ t,q>0 $ and $ x,y \in \bR \setminus \{ 0 \} $. 
Note that the second identity in \eqref{uq0xy} follows from Theorem \ref{thm: hitting}. 
Then $ p^0_t(x,y) $ is the continuous density of the transition probability 
for the killed process 
and $ u^0_q(x,y) $ is that of the resolvent.

\subsection{Excursion}

Since every point in $ \bR $ is regular for itself, 
the process admits the {\em local time} $ L(t,x) $ (see, e.g., \cite[Chap. V]{MR1406564}): 
$ L(t,x) $ is a measurable process such that 
\\ \quad 
{\rm (i)} 
for any $ x \in \bR $, $ t \mapsto L(t,x) $ is continuous $ P_0 $-almost surely 
(see, e.g., \cite[Proposition V.1.2]{MR1406564}. 
For joint continuity, see, e.g., \cite[Theorem V.3.15]{MR1406564}); 
\\ \quad 
{\rm (ii)} 
for $ P_0 $-almost all path, 
$ \displaystyle 
\int_0^t 1_A(X_s) \d s = \int_A L(t,x) \d x 
$ for all $ t \ge 0 $ and $ A \in \cB(\bR) $. 
\\
We simply write $ L_t = L(t,0) $. 
Denote its right-continuous inverse by $ \tau(l) $: 
\begin{align}
\tau(l) = \inf \{ t>0: L(t,0)>l \} . 
\n
\end{align}
Then the process $ (\tau(l):l \ge 0) $ is a subordinator 
such that 
\begin{align}
P_0[\e^{-q \tau(l)}] = \e^{-l/u_q(0)} 
, \qquad q>0 , \ l \ge 0 
\n
\end{align}
(see, e.g., \cite[pp. 131]{MR1406564}). 

Now we can apply It\^o's excursion theory 
(\cite{MR0402949}; see also \cite{MR1406564} and \cite{MR1138461} for details). 
We adopt the same notations as in \cite[Section 3]{YYY}. 
We denote its characteristic measure by $ \vn $ 
and call it the {\em excursion measure}. 
The measure $ \vn $ has its support 
on the set of {\em excursions away from the origin}: 
\begin{align}
\bE = \{ e \in \bD: 0 < \zeta(e) \le \infty \} 
\cap \{ e \in \bD : e(t) \neq 0 \ \text{for} \ 0<t<\zeta(e) \} . 
\n
\end{align}
We call an element $ e $ of $ \bE $ an {\em excursion path}. 
If $ \zeta(e)=\infty $, then such an excursion path $ e \in \bE $ is called 
a {\em final excursion}. 
Note that the measure $ \vn $ is $ \sigma $-finite on $ \cF_t $ for any $ t>0 $; 
in fact, $ \vn(\zeta>t)<\infty $ and $ \{ \zeta>t \} \in \cF_t $. 

\begin{theorem}[Markov property of {$ \vn $}] \label{thm: Markov}
For any $ t>0 $ 
and for any non-negative $ \cF_t $-measurable functional $ Z_t $ 
and any non-negative functional $ F $ on $ \bD $, 
it holds that 
\begin{align}
\vn [ Z_t F(X_{t+\cdot}) ] = \int \vn[Z_t;X_t \in \d x] P^0_x [F(X)] . 
\n
\end{align}
\end{theorem}

For the proof of Theorem \ref{thm: Markov}, 
see \cite[Theorem 6.3]{MR0402949} 
and also \cite[Theorem XII.4.1]{MR1725357}. 

The following theorem asserts that extremeness of the excursion measure 
is equivalent to germ triviality (see also It\^o \cite[Theorem 7.1]{MR0402949}). 

\begin{theorem} \label{thm: extremeness}
The excursion measure $ \vn $ is an extreme direction if and only if 
the germ $ \sigma $-field $ \cF_{0+} $ 
is trivial under $ \vn $, i.e., 
\begin{align}
\text{for any $ A \in \cF_{0+} $, it holds either that $ \vn(A)=0 $ or that $ \vn(A^c)=0 $.} 
\n
\end{align}
\end{theorem}

Although it is rather obvious, 
we give the proof for convenience of the reader. 

\begin{proof}
Suppose that $ \cF_{0+} $ is trivial under $ \vn $ 
and that $ \vm $ is absolutely continuous with respect to $ \vn $. 
We denote by $ D $ the Radon--Nikodym density on $ \cF_{\infty } $: 
$ \d \vm = D \d \vn $ on $ \cF_{\infty } $, 
and by $ D_t $ that on $ \cF_t $: 
$ \d \vm = D_t \d \vn $ on $ \cF_t $. 
Then, by the Markov property of $ \vm $ and $ \vn $, we have 
\begin{align}
\vm[Z_t F(X_{t+\cdot})] 
= \vm[Z_t P^0_{X_t}[F(X)]] 
= \vn[D_t Z_t P^0_{X_t}[F(X)]] 
= \vn[D_t Z_t F(X_{t+\cdot})] . 
\n
\end{align}
Hence we have $ D = D_t $ $ \vn $-almost everywhere, 
which implies that $ D $ is $ \cF_{0+} $-measurable $ \vn $-almost everywhere. 
Now we see that $ D $ is constant $ \vn $-almost everywhere. 
This proves that $ \vn $ is an extreme direction. 

Suppose that $ \cF_{0+} $ is not trivial under $ \vn $. 
Then there exists $ A \in \cF_{0+} $ such that $ \vn(A) \neq 0 $ 
and $ \vn(A^c) \neq 0 $. 
The measure $ \vn $ is decomposed into the sum 
\begin{align}
\d \vn = 1_A \d \vn + 1_{A^c} \d \vn, 
\n
\end{align}
where the measures $ 1_A \d \vn $ and $ 1_{A^c} \d \vn $ are mutually singular 
and they are elements of $ \cE $ (see \cite[Theorem 2]{MR859837}). 
Thus we see that $ \vn $ is not an extreme direction. 
The proof is complete. 
\qed
\end{proof}

\subsection{Duality} \label{sec: hitting}

Let $ \hat{P}_x $ denote the law of $ (-X_t) $ under $ P_{-x} $. 
Then the following duality holds: 
\begin{align}
\int f(x) P_x[g(X_t)] \d x = \int \hat{P}_x[f(X_t)] g(x) \d x 
\n
\end{align}
for any non-negative measurable functions $ f $ and $ g $. 
Since Theorem \ref{thm: hitting} is valid also for the dual process $ \{ (X_t),(\hat{P}_x) \} $, 
we have 
\begin{align}
\hat{P}_x \ebra{ \e^{-q T_{\{ 0 \}}} } = \frac{u_q(x)}{u_q(0)} 
, \qquad 
x \in \bR, \ q>0 . 
\label{hitting Lap trans dual}
\end{align}
The following theorem is due to 
Chen--Fukushima--Ying {\cite[Eq. (2.8)]{CFY}} 
and Fitzsimmons--Getoor {\cite[eq. (3.22)]{MR2247835}}, 
where the theorem has been proved in quite general settings. 

\begin{theorem} \label{lem: eq: master1}
Suppose that the conditions {\bf (A1)}-{\bf (A2)} are satisfied. 
Then, for any non-negative measurable function $ f $, it holds that 
\begin{align}
\int_0^{\infty } \e^{-qt} \vn[f(X_t)] \d t 
= \int f(x) \hat{P}_x \ebra{ \e^{-q T_{\{ 0 \}}} } \d x 
\label{eq: master1}
\end{align}
\end{theorem}

For the proof of Theorem \ref{lem: eq: master1}, 
see also \cite[Theorem 3.3]{YYY}. 

\begin{corollary}
Suppose that the conditions {\bf (A1)}-{\bf (A2)} are satisfied. Then 
\begin{align}
\int_0^{\infty } \e^{-qt} \vn(\zeta>t) \d t 
= \frac{1}{q u_q(0)} 
, \qquad q>0 . 
\label{eq: zeta>t}
\end{align}
\end{corollary}

\begin{proof}
Taking $ f(x) \equiv 1 $ in the identity \eqref{eq: master1}, 
we obtain the desired identity. 
\qed
\end{proof}

The following theorem relates the entrance law density 
with the hitting time density. 

\begin{theorem} \label{thm: rho}
Suppose that the conditions {\bf (A1)}-{\bf (A2)} are satisfied. Then 
there exists a bi-measurable function $ \rho(t,x) $ 
which is at the same time a space density of the entrance law 
\begin{align}
\vn(X_t \in \d x) = \rho(t,x) \d x 
\label{eq: ent law density}
\end{align}
and a time density of the law of the first hitting time for the dual process 
\begin{align}
\hat{P}_x(T_{\{ 0 \}} \in \d t;T_{\{ 0 \}}<\infty ) = \rho(t,x) \d t . 
\label{eq: first pas density}
\end{align}
That is, 
\begin{align}
\rho(t,x) = \frac{\vn(X_t \in \d x)}{\d x} 
= \frac{\hat{P}_x(T_{\{ 0 \}} \in \d t;T_{\{ 0 \}}<\infty )}{\d t} . 
\n
\end{align}
\end{theorem}

For the proof of Theorem \ref{thm: rho}, 
see \cite[Theorem 3.5]{YYY}.

\begin{corollary}
Suppose that the conditions {\bf (A1)}-{\bf (A2)} are satisfied. Then 
it holds that 
\begin{align}
p^0_t(x,y) = p_t(y-x) - \kbra{ p_{\cdot}(0)*\rho(\cdot,-x)*\rho(\cdot,y) }(t) 
\n
\end{align}
where $ f*g $ stands for the convolution of $ f $ and $ g $: $ f*g(t) = \int_0^t f(t-s) g(s) \d s $. 
\end{corollary}

\begin{proof}
The identity \eqref{uq0xy} can be written as 
\begin{align}
u^0_q(x,y) = u_q(y-x) - u_q(0) \cdot \frac{u_q(-x)}{u_q(0)} \cdot \frac{u_q(y)}{u_q(0)} . 
\n
\end{align}
Since we have 
\begin{align}
\int_0^{\infty } \e^{-qt} \rho(t,x) \d t = \frac{u_q(x)}{u_q(0)} 
, \qquad q>0 , \ x \in \bR , 
\n
\end{align}
we obtain the desired result. 
\qed
\end{proof}

\subsection{Continuous entrance property}

\begin{theorem}
Suppose that the conditions {\bf (A1)}-{\bf (A2)} are satisfied. Then 
\begin{align}
\vn(\{ X_0=0 \}^c) = 0 . 
\label{eq: conti ent}
\end{align}
\end{theorem}

We must be careful in this property \eqref{eq: conti ent}; 
in fact, it is not necessarily satisfied by general Markov processes with jumps, 
while it is obviously satisfied by diffusion processess. 

\begin{proof}
It is well-known (see, e.g., \cite[Theorem 6.31.5]{MR1739520}) that 
any function $ f $ of class $ C^2 $ with compact support 
belongs to the domain of the generator and 
\begin{align}
& \lim_{t \to 0+} \frac{P_x[f(X_t)]-f(x)}{t} 
\n \\
=& v f''(x) 
+ a f'(x) 
+ \int_{\bR} \kbra{ f(x+y)-f(x) - \frac{y}{1+y^2} f'(x) } \nu(\d y) 
\n
\end{align}
where $ v $, $ a $ and $ \nu $ have been introduced in \eqref{def theta}. 
Suppose that $ f $ is non-negative 
and its support is contained in $ \bR \setminus \{ 0 \} $. 
Then we have 
\begin{align}
\lim_{t \to 0+} \frac{1}{t} P_0[f(X_t)] 
= \int_{\bR} f(y) \nu(\d y) . 
\n
\end{align}
This implies that 
\begin{align}
C := \sup_{t>0} \frac{1}{t} P_0[f(X_t)] < \infty . 
\n
\end{align}
Now we have 
\begin{align}
q^2 \int u_q(x) f(x) \d x 
= q^2 \int_0^{\infty } \e^{-qt} P_0[f(X_t)] \d t 
\le C 
, \qquad q > 0 . 
\n
\end{align}
Thus, by \eqref{eq: master1}, we have 
\begin{align}
q \int_0^{\infty } \e^{-qt} \vn[f(X_t)] \d t 
= \frac{1}{q u_q(0)} \cdot q^2 \int u_q(x) f(x) \d x 
\le \frac{C}{q u_q(0)} , 
\n
\end{align}
which converges to 0 as $ q \to \infty $. 
This implies that 
\begin{align}
\liminf_{t \to 0+} \vn[f(X_t)] = 0 . 
\n
\end{align}
By Fatou's lemma, we obtain 
\begin{align}
\vn[f(X_0)] \le \liminf_{t \to 0+} \vn[f(X_t)] = 0 . 
\n
\end{align}
This proves $ \vn(\{ X_0=0 \}^c)=0 $. 
\qed
\end{proof}

\section{Preliminary facts: symmetric case} \label{sec: prel sym}

In what follows, we suppose that the following condition is satisfied: 
\begin{quote}
{\bf (A0)} 
The process is symmetric, 
i.e., $ P_x=\hat{P}_x $, $ x \in \bR $. 
\end{quote}
Then we see that 
\begin{align}
- \psi(\lambda) 
= \theta(\lambda) 
= v \lambda^2 + 2 \int_{(0,\infty )} \cbra{ 1-\cos \lambda x } \nu(\d x) 
, \qquad \lambda \in \bR 
\n
\end{align}
where the corresponding L\'evy measure is also symmetric: $ \nu(-\d x) = \nu(\d x) $.

\subsection{The Kesten--Bretagnolle condition in the symmetric case} \label{sec: KB}

\begin{lemma} \label{lem: cond A}
Suppose that the condition {\bf (A0)} is satisfied. Then 
the conditions {\bf (A1)}-{\bf (A2)} are satisfied if and only if 
\begin{align}
\int_0^{\infty } \frac{1}{1+\theta(\lambda)} \d \lambda < \infty , 
\label{eq: KBs}
\end{align}
and 
\begin{align}
\theta(\lambda) \to \infty 
\qquad \text{as} \qquad 
\lambda \to \infty . 
\label{eq: theta lambda conv to infty}
\end{align}
In this case, moreover, it holds that 
\begin{align}
\int_0^{\infty } \frac{\min \{ \lambda^2,1 \}}{\theta(\lambda)} \d \lambda < \infty . 
\label{eq: theta ibility}
\end{align}
\end{lemma}

\begin{proof}
(i) 
Suppose that the conditions {\bf (A1)}-{\bf (A2)} are satisfied. 
Then \eqref{eq: KB} of Theorem \ref{thm: Kesten} implies that \eqref{eq: KBs} holds. 
Note that 
\begin{align}
\e^{- \theta(\lambda)} 
= \int \e^{i \lambda x} p_1(x) \d x 
, \qquad \lambda \in \bR . 
\n
\end{align}
Since $ p_1(x) \ge 0 $ and $ \int p_1(x) \d x < \infty $, 
we may apply the Riemann--Lebesgue theorem to obtain 
$ \e^{- \theta(\lambda)} \to 0 $ as $ \lambda \to \infty $, 
which implies \eqref{eq: theta lambda conv to infty}. 

Suppose that $ \theta(\lambda_0)=0 $ for some $ \lambda_0 \in \bR \setminus \{ 0 \} $. 
Then we have $ P_0[\e^{i \lambda_0 X_t}] = 1 $, which implies that $ P_x(X_t=0)=1 $ 
for all $ t \ge 0 $. This contradicts the assumption {\bf (A2)}. 
Hence we obtain 
\begin{align}
\theta(\lambda) \neq 0 , 
\qquad \lambda \in \bR \setminus \{ 0 \} . 
\label{eq: not comp Poi}
\end{align}

It holds that 
\begin{align}
\frac{\theta(\lambda)}{\lambda^2} 
\ge& v + 2 \int_{(0,1)} \frac{1-\cos \lambda x}{(\lambda x)^2} x^2 \nu(\d x) 
\n \\
\to& 
v + \int_{(0,1)} x^2 \nu(\d x) 
\qquad \text{as $ \lambda \to 0 $} 
\label{eq: limit}
\end{align}
by the dominated convergence theorem. 
By \eqref{eq: KB2} of Theorem \ref{thm: Kesten}, we see that 
the limit in \eqref{eq: limit} is positive. 
Hence we see that there exists a positive constant $ C $ 
such that $ \theta(\lambda) \ge C \lambda^2 $ for small $ \lambda $. 
Combining it 
with \eqref{eq: KBs}, \eqref{eq: theta lambda conv to infty} and \eqref{eq: not comp Poi}, 
we obtain \eqref{eq: theta ibility}. 

(ii) 
Suppose that \eqref{eq: KBs} and \eqref{eq: theta lambda conv to infty} 
are satisfied. 
Suppose also that \eqref{eq: KB2} is {\em not} satisfied, 
i.e., that $ v=0 $ and $ \int_{(0,1)} x \nu(\d x) < \infty $. 
Then we have 
\begin{align}
\frac{\theta(\lambda)}{\lambda} 
\le& \int_{(0,1)} \frac{1-\cos \lambda x}{\lambda x} x \nu(\d x) 
+ \frac{4}{\lambda} \nu([1,\infty )) 
\n \\
\to& 0 
\qquad \text{as $ \lambda \to \infty $} 
\n
\end{align}
by the dominated convergence theorem. 
Hence we have $ \int_{\lambda_0}^{\infty } \frac{\d \lambda}{\theta(\lambda)} = \infty $ 
for any $ \lambda_0>0 $, 
which contradicts \eqref{eq: KBs} and \eqref{eq: theta lambda conv to infty}. 
The proof is now complete. 
\qed
\end{proof}

In what follows, we suppose that 
the condition {\bf (A)} is satisfied, 
i.e., that 
all of the conditions {\bf (A0)}-{\bf (A2)} are satisfied. 
Then we have 
\begin{align}
p_t(x) = \frac{1}{\pi} 
\int_0^{\infty } \cbra{ \cos \lambda x } \e^{-t \theta(\lambda)} \d \lambda 
\n
\end{align}
and 
\begin{align}
u_q(x) = \frac{1}{\pi} 
\int_0^{\infty } \frac{ \cos \lambda x }{q+\theta(\lambda)} \d \lambda . 
\n
\end{align}
Moreover, with the help of Lemma \ref{lem: cond A}, the function 
\begin{align}
h(x) 
= \lim_{q \to 0+} \kbra{ u_q(0)-u_q(x) } 
= \frac{1}{\pi} \int_0^{\infty } \frac{1-\cos \lambda x}{\theta(\lambda)} \d \lambda 
\label{eq: hx}
\end{align}
is well-defined.

\subsection{Recurrence and transience} \label{sec: rec and trans}

Recall that 
\begin{align}
\kappa 
= \lim_{q \to 0+} \frac{1}{u_q(0)} 
= \kbra{ \frac{1}{\pi} \int_0^{\infty } \frac{1}{\theta(\lambda)} \d \lambda }^{-1} 
\in [0,\infty ) . 
\n
\end{align}
We note that $ \int_1^{\infty } \frac{1}{\theta(\lambda)} \d \lambda < \infty $ 
by Lemma \ref{lem: cond A}. 
Now the following equivalence relations are well-known (see, e.g., \cite[Section I.4]{MR1406564}): 
\\ \quad (i) 
$ \kappa=0 $, or $ \int_0^1 \frac{1}{\theta(\lambda)} \d \lambda = \infty $, 
if and only if the process is {\em recurrent}, 
i.e., 
\begin{align}
P_x(T_{(-r,r)}<\infty )=1 
, \qquad x \in \bR , \ r>0 
\n
\end{align}
where $ T_{(-r,r)} $ stands for the first hitting time of the interval $ (-r,r) $; 
\\ \quad (ii) 
$ \kappa>0 $, or $ \int_0^1 \frac{1}{\theta(\lambda)} \d \lambda < \infty $, 
if and only if the process is {\em transient}, 
i.e., 
\begin{align}
P_x \cbra{ \lim_{t \to \infty }|X_t|=\infty } =1 
, \qquad x \in \bR. 
\n
\end{align}

\begin{theorem} \label{thm: transient}
Suppose that the condition {\bf (A)} is satisfied. Then 
\\ \quad {\rm (i)} 
If the process is recurrent, then 
\begin{align}
P_x(T_{\{ 0 \}} < \infty ) = 1 
, \qquad x \in \bR ; 
\n
\end{align}
\\ \quad {\rm (ii)} 
If the process is transient, then 
\begin{align}
P_x(T_{\{ 0 \}} = \infty ) = \kappa h(x) 
, \qquad x \in \bR . 
\n
\end{align}
\end{theorem}

\begin{proof}
On one hand, we have 
\begin{align}
\lim_{q \to 0+} P_x \ebra{ \e^{-q T_{\{ 0 \}}} } 
= P_x(T_{\{ 0 \}}<\infty ) . 
\n
\end{align}
On the other hand, we have 
\begin{align}
P_x \ebra{ \e^{-q T_{\{ 0 \}}} } 
= \frac{u_q(x)}{u_q(0)} 
= 1 - \frac{u_q(0)-u_q(x)}{u_q(0)} 
\stackrel{q \to 0+}{\longrightarrow} 1-\kappa h(x) . 
\n
\end{align}
This proves the claims (i) and (ii) at the same time. 
\qed
\end{proof}

\subsection{The distribution of the lifetime} \label{sec: dist of lifetime}

\begin{theorem} \label{thm: existence of lifetime law density}
Suppose that the condition {\bf (A)} is satisfied. Then 
there exists a completely monotone density $ \rho(t) $ such that 
\begin{align}
\vn(\zeta \in \d t) = \rho(t) \d t 
\qquad \text{on $ (0,\infty ) $} 
\label{eq: lifetime-density}
\end{align}
and that 
\begin{align}
\vn(\zeta=\infty ) = \kappa 
\label{eq: final excursion}
\end{align}
\end{theorem}

\begin{proof}
By \eqref{eq: zeta>t}, we obtain 
\begin{align}
\vn(\zeta=\infty ) 
= \lim_{q \to 0+} \frac{1}{u_q(0)} = \kappa . 
\n
\end{align}

Let us introduce the following positive Borel measure $ \sigma $ on $ [0,\infty ) $: 
\begin{align}
\sigma(A) = \frac{1}{\pi} \int_0^{\infty } 1_A(\theta(\lambda)) \d \lambda 
, \qquad A \in \cB([0,\infty )) . 
\n
\end{align}
Then we have 
\begin{align}
\int_{[0,\infty )} \frac{1}{1+\xi} \sigma(\d \xi) 
= \frac{1}{\pi} \int_0^{\infty } \frac{1}{1+\theta(\lambda)} \d \lambda 
< \infty . 
\n
\end{align}
In particular, we see that $ \sigma $ is a Radon measure. 
Since we have 
\begin{align}
u_q(0) = \int_{[0,\infty )} \frac{1}{q + \xi} \sigma(\d \xi) 
, \qquad q>0 , 
\n
\end{align}
it is well-known (see, e.g., \cite[Chapter II]{MR0486556}) 
that there exists another Radon measure $ \sigma^* $ on $ [0,\infty ) $ 
such that 
\begin{align}
\int_{[0,\infty )} \frac{1}{1+\xi} \sigma^*(\d \xi) < \infty 
\n
\end{align}
and 
\begin{align}
\frac{1}{q u_q(0)} = \int_{[0,\infty )} \frac{1}{q + \xi} \sigma^*(\d \xi) 
, \qquad q>0 . 
\label{1/quq0}
\end{align}
Combining \eqref{eq: zeta>t} and \eqref{1/quq0}, we obtain 
\begin{align}
\vn(\zeta>t) = \int_{[0,\infty )} \e^{-t \xi} \sigma^*(\d \xi) 
, \qquad t>0 . 
\label{eq: zeta and sigma*}
\end{align}
Therefore we conclude that 
$ \vn(\zeta \in \d t) = \rho(t) \d t $ where 
\begin{align}
\rho(t) = \int_{(0,\infty )} \e^{-t \xi} \xi \sigma^*(\d \xi) 
, \qquad t>0. 
\label{rhot}
\end{align}
It is obvious by definition that 
the function $ \rho(t) $ is completely monotone. 
Now the proof is complete. 
\qed
\end{proof}

\subsection{Another expression of $ \rho(t,x) $}

The following theorem gives another expression of the density $ \rho(t,x) $: 

\begin{theorem} \label{thm: another}
Suppose that the condition {\bf (A)} is satisfied. Then 
\begin{align}
\rho(t,x) = \frac{1}{t} 
\int_0^t \kbra{ R(t-s) p_s(x) + (t-s)R'(t-s) p_s(x) + R(t-s) s \frac{\d}{\d s} p_s(x) } \d s 
\n
\end{align}
where $ R(t) = \vn(\zeta>t) = \int_{[0,\infty )} \e^{-t\xi} \sigma^*(\d \xi) $. 
\end{theorem}

Theorem \ref{thm: another} may be proved 
in the same way as \cite[Proposition 4.2]{YY}, 
so we omit the proof.

\section{Harmonic function for the killed process} \label{sec: h-func}

Set 
\begin{align}
h_q(x) = u_q(0)-u_q(x) 
= \frac{1}{\pi} \int_0^{\infty } \frac{1-\cos \lambda x}{q + \theta(\lambda)} \d \lambda 
, \qquad q>0 , \ x \in \bR . 
\label{hqx}
\end{align}
We note that $ h_q(x) \ge 0 $ 
and that, for each $ x \in \bR $, 
the function $ q \mapsto h_q(x) $ increases as $ q>0 $ decreases. 
Recall that 
\begin{align}
h(x) 
= \lim_{q \to 0+} h_q(x) 
= \frac{1}{\pi} \int_0^{\infty } \frac{1-\cos \lambda x}{\theta(\lambda)} \d \lambda 
, \qquad x \in \bR . 
\n
\end{align}

\subsection{Harmonicity of $ h(x) $}

To prove that the function $ h(x) $ is harmonic, 
we need the following 

\begin{lemma} \label{lem: uq and h}
Suppose that the condition {\bf (A)} is satisfied. Then 
it holds that 
\begin{align}
\int u_q(y-x) h(y) \d y = \frac{h(x) + u_q(x)}{q} 
, \qquad q>0 , \ x \in \bR , 
\label{eq: uq and h}
\end{align}
and, consequently, that 
\begin{align}
\int p_t(y-x) h(y) \d y = h(x) + \int_0^t p_s(x) \d s 
, \qquad t>0 , \ x \in \bR . 
\label{eq: pt and h}
\end{align}
\end{lemma}

\begin{proof}
Note that 
\begin{align}
\int u_q(y-x) h(y) \d y 
= \lim_{r \to 0+} 
\int u_q(y-x) h_r(y) \d y 
\n
\end{align}
by the monotone convergence theorem. 
Let $ r $ be such that $ 0<r<q $. Then 
\begin{align}
\int u_q(y-x) h_r(y) \d y 
=& \int u_q(y-x) \{ u_r(0)-u_r(y) \} \d y 
\n \\
=& \frac{u_r(0)}{q} - \int u_q(y-x) u_r(-y) \d y 
\n
\end{align}
where we used the symmetry: $ u_r(y)=u_r(-y) $. 
By the resolvent equation \eqref{eq: resolv eq} with $ z=0 $, 
the last equation becomes 
\begin{align}
\frac{u_r(0)}{q} - \frac{1}{q-r} \kbra{ u_r(x) - u_q(x) } 
= \frac{h_r(x)}{q} - \frac{ru_r(x)}{q(q-r)} + \frac{u_q(x)}{q-r} . 
\n
\end{align}
Letting $ r \to 0+ $, 
we have $ h_r(x) \to h(x) $ and $ ru_r(x) \to 0 $, 
and hence we see that 
the last equation tends to the right hand side of \eqref{eq: uq and h}. 
Now \eqref{eq: pt and h} is obvious, 
and the proof is complete. 
\qed
\end{proof}

Now we proceed to the proof of Theorem \ref{thm: harmonic}. 

%\begin{proof}[Proof of Theorem \ref{thm: harmonic}]
\noindent 
{\it Proof of Theorem \ref{thm: harmonic}} 
In order to prove the assertion, 
it suffices to show that 
\begin{align}
q U^0_q h = h 
, \qquad q>0 . 
\n
\end{align}
By definition of $ U^0_q $, we have 
\begin{align}
U^0_q h(x) 
=& \int u^0_q(x,y) h(y) \d y 
\n \\
=& \int u_q(y-x) h(y) \d y - \frac{u_q(x)}{u_q(0)} \int u_q(y) h(y) \d y . 
\n
\end{align}
By Lemma \ref{lem: uq and h}, we obtain 
\begin{align}
U^0_q h(x) 
=& \frac{h(x) + u_q(x)}{q} - \frac{u_q(x)}{u_q(0)} \cdot \frac{u_q(0)}{q} 
= \frac{h(x)}{q} . 
\n
\end{align}
This completes the proof. 
\qed
%\end{proof}

Now let us prove Theorem \ref{thm: vn hXt = 1}. 

%\begin{proof}[Proof of Theorem \ref{thm: vn hXt = 1}]
\noindent 
{\it Proof of Theorem \ref{thm: vn hXt = 1}} 
By Lemma \ref{lem: uq and h}, we have 
\begin{align}
\int h(y) \frac{u_q(y)}{u_q(0)} \d y = \frac{1}{q} 
, \qquad q>0 . 
\label{eq: to show vn hXt = 1 -2}
\end{align}
By \eqref{eq: ent law density}, \eqref{eq: first pas density} and \eqref{hitting Lap trans}, 
we see that the identity \eqref{eq: to show vn hXt = 1 -2} is equivalent to 
\begin{align}
\int_0^{\infty } \e^{-qt} \vn [h(X_t)] \d t = \frac{1}{q} 
, \qquad q>0 . 
\label{eq: to show vn hXt = 1}
\end{align}
Therefore we obtain the desired identity. 
\qed
%\end{proof}

\subsection{Several properties of $ h(x) $}

Let us study several properties of $ h(x) $. 

\begin{lemma} \label{lem: properties of hx}
Suppose that the condition {\bf (A)} is satisfied. Then 
the following assertions hold: 
\\ \quad {\rm (i)} 
$ h(x) $ is continuous; 
\\ \quad {\rm (ii)} 
$ h(0)=0 $; 
\\ \quad {\rm (iii)} 
$ h(x)>0 $ for all $ x \in \bR \setminus \{ 0 \} $; 
\\ \quad {\rm (iv)} 
$ \displaystyle \lim_{|x| \to \infty } h(x) 
= \frac{1}{\kappa} 
= \frac{1}{\pi} \int_0^{\infty } \frac{1}{\theta(\lambda)} \d \lambda \in (0,\infty ] $. 
\end{lemma}

\begin{proof}
The assertion (i) is obvious by the dominated convergence theorem. 
The assertions (ii) and (iii) are obvious by definition. 
Let us prove the assertion (iv). 

\underline{Transient case}. 
Since $ \int_0^{\infty } \frac{1}{\theta(\lambda)} \d \lambda < \infty $, 
we may apply the Riemann--Lebesgue theorem to obtain 
\begin{align}
\lim_{|x| \to \infty } \int_0^{\infty } \frac{\cos \lambda x}{\theta(\lambda)} \d \lambda 
= 0 . 
\n
\end{align}
This proves that $ \lim_{|x| \to \infty } h(x) = 1/\kappa $. 

\underline{Recurrent case}. 
We have $ \int_0^{\infty } \frac{1}{\theta(\lambda)} \d \lambda = \infty $. 
Let $ \eps>0 $. 
Since $ \int_{\eps}^{\infty } \frac{1}{\theta(\lambda)} \d \lambda < \infty $, 
we may apply the Riemann--Lebesgue theorem to obtain 
\begin{align}
\lim_{|x| \to \infty } \int_{\eps}^{\infty } \frac{\cos \lambda x}{\theta(\lambda)} \d \lambda 
= 0 . 
\n
\end{align}
Hence we obtain 
\begin{align}
\liminf_{|x| \to \infty } \pi h(x) 
\ge \int_{\eps}^{\infty } \frac{1}{\theta(\lambda)} \d \lambda . 
\n
\end{align}
Letting $ \eps \to 0+ $, we obtain $ \lim_{|x| \to \infty } h(x) = \infty $. 
\qed
\end{proof}

\begin{lemma} \label{lem: bdd}
Suppose that the condition {\bf (A)} is satisfied. Then 
\begin{align}
\sup_{|x|<1} \frac{|x|}{h(x)} < \infty . 
\n
\end{align}
\end{lemma}

\begin{proof}
By \eqref{def theta}, we have 
\begin{align}
\theta(\lambda) \le v \lambda^2 + \lambda^2 \int_{(0,1)} x^2 \nu(\d x) + 4 \nu([1,\infty )) 
, \qquad \lambda \in \bR . 
\n
\end{align}
Now we see that there exists a constant $ C $ such that 
\begin{align}
\theta(\lambda) \le C \lambda^2 
, \qquad \lambda>1 . 
\n
\end{align}
Hence we obtain 
\begin{align}
\pi h(x) 
\ge 
\int_1^{\infty } \frac{1-\cos \lambda x}{\theta(\lambda)} \d \lambda 
\ge \frac{1}{C} 
\int_1^{\infty } \frac{1-\cos \lambda x}{\lambda^2} \d \lambda 
= \frac{|x|}{C} 
\int_{|x|}^{\infty } \frac{1-\cos \lambda }{\lambda^2} \d \lambda . 
\n
\end{align}
For $ |x|<1 $, we obtain 
\begin{align}
\pi h(x) 
\ge \frac{|x|}{C} 
\int_1^{\infty } \frac{1-\cos \lambda }{\lambda^2} \d \lambda . 
\n
\end{align}
Since $ \int_1^{\infty } \frac{1-\cos \lambda }{\lambda^2} \d \lambda >0 $, 
we complete the proof. 
\qed
\end{proof}

\begin{lemma} \label{lem: conv to 0}
Suppose that the condition {\bf (A)} is satisfied. 
Then, for any fixed $ q>0 $, 
\begin{align}
\lim_{x \to 0} \frac{ h(x) - h_q(x) }{h(x)} = 0 . 
\label{h-trans:eq1}
\end{align}
\end{lemma}

\begin{proof}
By definitions, we have 
\begin{align}
h(x) - h_q(x) 
= \frac{1}{\pi} \int_0^{\infty } 
\frac{1-\cos \lambda x}{\theta(\lambda)} 
\frac{q}{q+\theta(\lambda)} 
\d \lambda 
\ge 0. 
\label{hx-hqx}
\end{align}
Let $ \eps>0 $ be arbitrary. 
Since $ \theta(\lambda) \to \infty $ as $ \lambda \to \infty $ (Lemma \ref{lem: cond A}), 
there exists a constant $ L>0 $ such that $ \frac{q}{q+\theta(\lambda)}<\eps $ 
for all $ \lambda>L $. 
Now we have 
\begin{align}
h(x) - h_q(x) 
\le& \frac{1}{\pi} \int_0^L \frac{1-\cos \lambda x}{\theta(\lambda)} \d \lambda 
+ \frac{\eps}{\pi} \int_L^{\infty } \frac{1-\cos \lambda x}{\theta(\lambda)} \d \lambda 
\n \\
\le& \frac{x^2}{2 \pi} \int_0^L \frac{\lambda^2}{\theta(\lambda)} \d \lambda 
+ \eps h(x) 
\n
\end{align}
for all $ x \in \bR $. 
Hence, by Lemma \ref{lem: bdd}, we obtain 
\begin{align}
\limsup_{x \to 0} \frac{h(x) - h_q(x)}{h(x)} \le \eps . 
\n
\end{align}
Since $ \eps>0 $ is arbitrary, 
we obtain \eqref{h-trans:eq1}. 
\qed
\end{proof}

\subsection{Green function for the killed process}

The function 
\begin{align}
u^0_0(x,y) 
= \lim_{q \to 0+} u^0_q(x,y) 
, \qquad (x,y) \in \bR^2 \setminus \{ (0,0) \} 
\n
\end{align}
is called the {\em Green function for the killed process}. 
The following formula can also be found 
in Eisenbaum--Kaspi--Marcus--Rosen--Shi \cite[Theorem 6.1]{MR1813843}: 

\begin{lemma} \label{lem: Green for killed}
Suppose that the condition {\bf (A)} is satisfied. Then 
\begin{align}
u^0_0(x,y) 
= 
h(x)+h(y)-h(y-x) - \kappa h(x) h(y) 
, \qquad x,y \in \bR \setminus \{ 0 \} . 
\label{eq: Green func for killed process}
\end{align}
\end{lemma}

\begin{proof}
Noting that 
\begin{align}
u^0_q(x,y) 
= h_q(x) + h_q(y) - h_q(y-x) - \frac{h_q(x)h_q(y)}{u_q(0)} , 
\label{eq: u harm qxy}
\end{align}
we obtain the desired result by letting $ q \to 0+ $. 
\qed
\end{proof}

\section{The $ h $-path process of the killed process} \label{sec: h-trans}

Recall that the $ h $-path process of the killed process has been introduced as 
\begin{align}
\d P^{\harm}_x |_{\cF_t} = 
\begin{cases}
\displaystyle 
\frac{h(X_t)}{h(x)} \d P^0_x |_{\cF_t} 
, \quad & x \in \bR \setminus \{ 0 \} , 
\\
h(X_t) \d \vn |_{\cF_t} 
, \quad & x=0 . 
\end{cases}
\n
\end{align}
It is obvious by definition that 
$ P^{\harm}_x(T_{\{ 0 \}}=+\infty ) = 1 $ for any $ x \in \bR $. 
For $ t>0 $, we define 
\begin{align}
p^{\harm}_t(x,y) = 
\frac{p^0_t(x,y)}{h(x) h(y)} 
, \qquad x,y \in \bR \setminus \{ 0 \} 
\n
\end{align}
and 
\begin{align}
p^{\harm}_t(x,0) = 
p^{\harm}_t(0,x) = 
\frac{\rho(t,x)}{h(x)} 
, \qquad x \in \bR \setminus \{ 0 \} . 
\n
\end{align}
Then $ p^{\harm}_t(x,y) $ is a density 
of the transition probability 
of the process $ \{ (X_t:t \ge 0),(P^{\harm}_x:x \in \bR) \} $ 
with respect to the symmetrizing measure $ h(y)^2 \d y $: 
\begin{align}
P^{\harm}_x(X_t \in A) = \int_A p^{\harm}_t(x,y) h(y)^2 \d y 
\n
\end{align}
for $ t>0 $, $ x \in \bR $ and $ A \in \cB(\bR) $. 
For $ q \ge 0 $, we define 
\begin{align}
u^{\harm}_q(x,y) = 
\frac{u^0_q(x,y)}{h(x) h(y)} 
, \qquad x,y \in \bR \setminus \{ 0 \} 
\n
\end{align}
and 
\begin{align}
u^{\harm}_q(x,0) = 
u^{\harm}_q(0,x) = 
\frac{u_q(x)}{h(x) u_q(0)} 
, \qquad x \in \bR \setminus \{ 0 \} . 
\n
\end{align}
Then $ u^{\harm}_q(x,y) $ is a density 
of the resolvent kernel: 
\begin{align}
P^{\harm}_x \ebra{ \int_0^{\infty } \e^{-qt} 1_A(X_t) \d t } 
= \int_A u^{\harm}_q(x,y) h(y)^2 \d y 
\n
\end{align}
for $ q \ge 0 $, $ x \in \bR $ and $ A \in \cB(\bR) $.

\subsection{Chapman--Kolmogorov identities}

We define 
\begin{align}
p^{\harm}_t(0,0) = \rho(t) 
, \qquad t>0 
\n
\end{align}
where the function $ \rho(t) $ has been introduced in \eqref{rhot}. 
Then the formulae \eqref{eq: zeta>t} and \eqref{eq: final excursion} imply that 
\begin{align}
\int_0^{\infty } \cbra{ 1-\e^{-qt} } p^{\harm}_t(0,0) \d t 
= \frac{1}{u_q(0)} - \kappa 
, \qquad q>0 . 
\label{LT of p never t00}
\end{align}

\begin{proposition} \label{prop: CK}
Suppose that the condition {\bf (A)} is satisfied. Then 
the following Chapman--Kolmogorov identities hold: 
\begin{align}
\int p^{\harm}_s(x,y) p^{\harm}_t(y,z) h(y)^2 \d y = p^{\harm}_{s+t}(x,z) 
, \qquad s,t>0 , \ x,z \in \bR . 
\label{CK}
\end{align}
\end{proposition}

\begin{proof}%[Proof of Proposition \ref{prop: CK}]
The identity \eqref{CK} is immediate 
in the case where $ x,z \in \bR \setminus \{ 0 \} $ by the Markov property of the killed process, 
and in the case where $ xz=0 $ except where $ x=z=0 $ by the Markov property of the excursion measure. 
Hence it suffices to prove the identity \eqref{CK} in the case where $ x=z=0 $. 

Let $ q,r>0 $ with $ q \neq r $. 
On the one hand, letting $ x=z=0 $ in the resolvent equation \eqref{eq: resolv eq} 
and using the symmetry, we have 
\begin{align}
\int u_q(y) u_r(y) \d y = \frac{1}{q-r} \kbra{ u_r(0) - u_q(0) } . 
\n
\end{align}
Dividing both sides by $ u_q(0) u_r(0) $, we obtain 
\begin{align}
\int u^{\harm}_q(0,y) u^{\harm}_r(y,0) h(y)^2 \d y 
= \frac{1}{q-r} \kbra{ \frac{1}{u_q(0)} - \frac{1}{u_r(0)} } . 
\label{DLT1}
\end{align}
On the other hand, 
we compute the double Laplace transform of $ p^{\harm}_{s+t}(0,0) $ as 
\begin{align}
& \int_0^{\infty } \d s \int_0^{\infty } \d t \e^{-qs-rt} p^{\harm}_{s+t}(0,0) 
\n \\
=& \frac{1}{q-r} \int_0^{\infty } \cbra{ \e^{-ru} - \e^{-qu} } p^{\harm}_u(0,0) \d u 
\n \\
=& \frac{1}{q-r} \kbra{ \int_0^{\infty } \cbra{ 1 - \e^{-qu} } p^{\harm}_u(0,0) \d u 
- \int_0^{\infty } \cbra{ 1 - \e^{-ru} } p^{\harm}_u(0,0) \d u } 
\n \\
=& \frac{1}{q-r} \kbra{ \frac{1}{u_q(0)} - \frac{1}{u_r(0)} } . 
\label{DLT2}
\end{align}
where we used the identity \eqref{LT of p never t00}. 
This shows that the right hand side of \eqref{DLT1} 
is the double Laplace transform of $ p^{\harm}_{s+t}(0,0) $. 
Since the left hand side of \eqref{DLT1} 
is the double Laplace transform of $ \int p^{\harm}_s(0,y) p^{\harm}_t(y,0) h(y)^2 \d y $, 
we obtain the desired identity. 
\qed
\end{proof}

\subsection{Lifetime disintegration formula} \label{sec: proof of Lifetime disinteg}

Let us define the law of the {\em bridge process} 
\begin{align}
P^{\harm}_{0,0;t}(\cdot) 
:= 
P^{\harm}_0(\cdot | X_{t-}=0) 
, \qquad t > 0 
\n
\end{align}
as the unique probability measure carried on the set 
$ \{ X_{t-}=0 , \, \zeta = t \} $ 
such that 
\begin{align}
\left. \d P^{\harm}_{0,0;t} \right|_{\cF_s} 
= \frac{p^{\harm}_{t-s}(X_s,0)}{p^{\harm}_t(0,0)} \left. \d P^{\harm}_0 \right|_{\cF_s} 
, \qquad s \in (0,t) . 
\n
\end{align}
Then the process $ \{ (X_s:0 \le s \le t),P^{\harm}_{0,0;t} \} $ 
is a time-inhomogeneous Markov process 
with its transition probability given by 
\begin{align}
P^{\harm}_{0,0;t}(X_v \in \d b | X_u=a) 
= \frac{p^{\harm}_{v-u}(a,b) p^{\harm}_{t-v}(b,0)}{p^{\harm}_{t-u}(a,0)} h(b)^2 \d b 
, \qquad 0<u<v<t , \ a,b \in \bR . 
\n
\end{align}
We call the process the {\em bridge} of the process $ \{ (X_s:s \ge 0),P^{\harm}_0 \} $ 
given $ X_{t-}=0 $. 
For further properties of the bridge process, see Fitzsimmons--Pitman--Yor \cite{FPY}.

Now we are in a position to prove Theorem \ref{thm: disinteg by lifetime}. 

%\begin{proof}[Proof of Theorem \ref{thm: disinteg by lifetime}]
\noindent 
{\it Proof of Theorem \ref{thm: disinteg by lifetime}} 
We consider a cylinder set $ \Gamma $ of $ \bE $ of the form 
\begin{align}
\Gamma = \kbra{ X_{t_1} \in A_1, \ldots, X_{t_n} \in A_n } 
\label{cylinder2}
\end{align}
for some sequence $ 0<t_1<t_2<\cdots<t_n<\infty $ 
and some sets $ A_1,\ldots,A_n \in \cB(\bR \setminus \{ 0 \}) $. 

Let $ t>t_n $. 
By the Markov property of the excursion measure, we have 
\begin{align}
\vn(\Gamma \cap \{ t<\zeta<\infty \}) 
= \vn \ebra{ 1_{\Gamma } P_{X_{t_n}}(t-t_n<T_{\{ 0 \}}<\infty ) } . 
\n
\end{align}
By the definition of $ P^{\harm}_0 $, we have 
\begin{align}
\vn(\Gamma \cap \{ t<\zeta<\infty \}) 
= P^{\harm}_0 \ebra{ 1_{\Gamma} \frac{P_{X_{t_n}}(t-t_n<T_{\{ 0 \}}<\infty )}{h(X_{t_n})} } . 
\n
\end{align}
Note that, by Theorem \ref{thm: rho}, we have 
\begin{align}
P_y(t-t_n<T_{\{ 0 \}}<\infty ) 
= \int_t^{\infty } \rho(s-t_n,y) \d s . 
\n
\end{align}
Hence we obtain 
\begin{align}
\vn(\Gamma \cap \{ t<\zeta<\infty \}) 
=& \int_t^{\infty } \d s P^{\harm}_0 \ebra{ 1_{\Gamma} \frac{\rho(s-t_n,X_{t_n})}{h(X_{t_n})} } 
\n \\
=& \int_t^{\infty } \d s 
P^{\harm}_0 \ebra{ 1_{\Gamma} p^{\harm}_{s-t_n}(X_{t_n},0) } 
\n \\
=& \int_t^{\infty } \d s p^{\harm}_s(0,0) 
P^{\harm}_0 \ebra{ 1_{\Gamma} \frac{p^{\harm}_{s-t_n}(X_{t_n},0)}{p^{\harm}_s(0,0)} } 
\n \\
=& \int_t^{\infty } \d s p^{\harm}_s(0,0) P^{\harm}_0 (\Gamma|X_{s-}=0) . 
\n
\end{align}
Letting $ t \to t_n+ $, we have 
\begin{align}
\vn(\Gamma \cap \{ \zeta<\infty \}) 
= \int_0^{\infty } 
P^{\harm}_0 \cbra{ \Gamma | X_{s-}=0 } p^{\harm}_s(0,0) \d s , 
\label{eq: zeta<infty}
\end{align}
since $ \Gamma \subset \{ \zeta>t_n \} $ 
and $ P^{\harm}_0 ( \zeta=s | X_{s-}=0 ) = 1 $. 

Note that, by Theorem \ref{thm: transient}, we have 
\begin{align}
P_y(T_{\{ 0 \}}= \infty ) = \kappa h(y) 
\n
\end{align}
whichever the process is recurrent or transient. 
Hence we obtain 
\begin{align}
\vn(\Gamma \cap \{ \zeta=\infty \}) 
= \kappa P^{\harm}_0(\Gamma) . 
\label{eq: zeta=infty}
\end{align}

Therefore, combining \eqref{eq: zeta<infty} and \eqref{eq: zeta=infty}, 
we obtain 
\begin{align}
\vn(\Gamma) = 
\int_0^{\infty } 
P^{\harm}_0 \cbra{ \Gamma | X_{s-}=0 } p^{\harm}_s(0,0) \d s 
+ 
\kappa P^{\harm}_0(\Gamma) . 
\n
\end{align}
Noting that the set of all cylinder sets $ \Gamma $ of the form \eqref{cylinder2} 
generate the whole $ \sigma $-field, 
we obtain the desired result. 
\qed
%\end{proof}

\subsection{Green function for the $ h $-path process}

The function 
\begin{align}
u^{\harm}_0(x,y) 
= \lim_{q \to 0+} u^{\harm}_q(x,y) 
, \qquad (x,y) \in \bR^2 \setminus \{ (0,0) \} 
\n
\end{align}
is called the {\em Green function for the $ h $-path process}. 

\begin{lemma}
Suppose that the condition {\bf (A)} is satisfied. Then 
\begin{align}
u^{\harm}_0(x,y) 
= 
\frac{1}{h(x)h(y)} \kbra{ h(x)+h(y)-h(y-x) } - \kappa 
, \qquad x,y \in \bR \setminus \{ 0 \} , 
\n
\end{align}
and that 
\begin{align}
u^{\harm}_0(x,0) = u^{\harm}_0(0,x) 
= 
\frac{1}{h(x)} - \kappa 
, \qquad x \in \bR \setminus \{ 0 \} . 
\n
\end{align}
\end{lemma}

\begin{proof}
The first identity is obvious from Lemma \ref{lem: Green for killed}. 
If $ x \in \bR \setminus \{ 0 \} $, 
then we have 
\begin{align}
h(x) u^{\harm}_q(x,0) 
= 
\frac{ u_q(0)-h_q(x) }{u_q(0)} , 
\label{eq: u harm qx}
\end{align}
and hence we obtain the second identity. 
\qed
\end{proof}

\section{Key lemmas} \label{sec: key lemmas}

We need several lemmas for later use.

\subsection{Regularity of $ h_q(x) $}

\begin{lemma}
Suppose that the conditions {\bf (A)} and {\bf (T)} are satisfied. 
For $ q \ge 0 $, set 
\begin{align}
\Theta_q(\d \eta) = \frac{1}{\pi(q+\theta(\eta))^2} \d \theta(\eta) 
\qquad \text{on $ (\lambda_0,\infty ) $} . 
\n
\end{align}
Then 
\begin{align}
\int_{(\lambda_0,\infty )} \eta \Theta_q(\d \eta) < \infty 
, \qquad q \ge 0 . 
\n
\end{align}
\end{lemma}

\begin{proof}
Since 
the function $ \lambda \mapsto \theta(\lambda) $ increases as $ \lambda>0 $ increases 
(the assumption {\bf (T)}) 
and 
since $ \theta(\lambda) \to \infty $ as $ \lambda \to \infty $ (Lemma \ref{lem: cond A}), 
we see that 
$ \Theta_q(\d \eta) $ is well-defined as a positive Borel measure on $ (\lambda_0,\infty ) $. 
Since 
\begin{align}
\int_{(\lambda_0,\infty )} \eta \Theta_q(\d \eta) 
=& \int_{(\lambda_0,\infty )} \Theta_q(\d \eta) 
\kbra{ \lambda_0 + \int_{\lambda_0}^{\eta} \d \lambda } 
\n \\
=& 
\frac{\lambda_0}{\pi(q+\theta(\lambda_0))} 
+ \int_{\lambda_0}^{\infty } \d \lambda \int_{(\lambda,\infty )} \Theta_q(\d \eta) 
\n \\
=& 
\frac{\lambda_0}{\pi(q+\theta(\lambda_0))} 
+ \frac{1}{\pi} \int_{\lambda_0}^{\infty } \frac{1}{q+\theta(\lambda)} \d \lambda 
< \infty , 
\label{ibility of F}
\end{align}
we complete the proof. 
\qed
\end{proof}

Set $ h_0 = h $. Then we have 
\begin{align}
h_q(x) = \frac{1}{\pi} \int_0^{\infty } \frac{1-\cos \lambda x}{q+\theta(\lambda)} \d \lambda 
, \qquad q \ge 0, \ x \in \bR . 
\n
\end{align}

\begin{lemma} \label{lem: estimate}
Suppose that the conditions {\bf (A)} and {\bf (T)} are satisfied. 
Let $ q \ge 0 $ be fixed. 
Then the following statements hold: 
\\ \quad 
{\rm (i)} 
There exists a constant $ C_q $ such that, for any $ x,y $ such that $ 0<2|x|<|y| $, 
\begin{align}
\abra{ \frac{h_q(y-x) - h_q(y)}{x} } \le C_q \cbra{ |y| + \frac{1}{|y|} } ; 
\label{eq: ineq Cq y-x}
\end{align}
\\ \quad 
{\rm (ii)} 
For any $ y \in \bR \setminus \{ 0 \} $, it holds that 
\begin{align}
\lim_{x \to 0} \frac{h_q(y+x) + h_q(y-x) - 2 h_q(y)}{x} = 0 ; 
\label{eq: limit = 0}
\end{align}
\\ \quad 
{\rm (iii)} 
For any $ \eps>0 $, it holds that 
\begin{align}
\lim_{x \to 0} \frac{1}{x} 
\int_{-\eps}^{\eps} \kbra{h_q(y+x) + h_q(y-x) - 2 h_q(y)} \d y = 0 . 
\label{eq: limit int = 0}
\end{align}
\end{lemma}

\begin{proof}%[Proof of Lemma \ref{lem: estimate}]
(i) 
For $ v \in \bR $, we split $ h_q(v) $ into the sum of 
\begin{align}
h^{(1)}_q(v) = \frac{1}{\pi} 
\int_0^{\lambda_0} \frac{1-\cos \lambda v}{q+\theta(\lambda)} \d \lambda 
\qquad \text{and} \qquad 
h^{(2)}_q(v) = \frac{1}{\pi} 
\int_{\lambda_0}^{\infty } \frac{1-\cos \lambda v}{q+\theta(\lambda)} \d \lambda . 
\n
\end{align}

We note that 
\begin{align}
\frac{h^{(1)}_q(y-x) - h^{(1)}_q(y)}{x} 
= \int_0^{\lambda_0} A^{(1)}_{\lambda,y}(x) 
\frac{\lambda^2}{q+\theta(\lambda)} \d \lambda 
\n
\end{align}
where 
\begin{align}
A^{(1)}_{\lambda,y}(x) 
= \frac{\cos(\lambda y) - \cos (\lambda(y-x)) }{\pi \lambda^2 x} . 
\n
\end{align}
Since we have 
\begin{align}
\abra{ A^{(1)}_{\lambda,y}(x) } 
= \frac{1}{\pi \lambda^2 |x|} \abra{ \int_0^x (-\lambda) \sin (\lambda(y-v)) \d v } 
\le \frac{2|y|}{\pi} , 
\label{eq: ineq A(1)}
\end{align}
we obtain 
\begin{align}
\abra{ \frac{h^{(1)}_q(y-x) - h^{(1)}_q(y)}{x} } 
\le \frac{2|y|}{\pi} \int_0^{\lambda_0} \frac{\lambda^2}{\theta(\lambda)} \d \lambda . 
\n
\end{align}

Define a function $ \varphi $ as $ \varphi(0)=1 $ and 
\begin{align}
\varphi(v) = \frac{\sin v}{v} 
, \qquad v \neq 0 . 
\n
\end{align}
Since $ \varphi'(v) = \frac{\cos v}{v} - \frac{\sin v}{v^2} $ 
and $ |\varphi'(v)| \le 2/|v| $ for all $ v \neq 0 $, 
we have 
\begin{align}
\abra{ \frac{\varphi(\lambda (y-x)) - \varphi(\lambda y)}{x} } 
= \abra{ \frac{1}{x} \int_0^x (-\lambda) \varphi'(\lambda (y-v)) \d v } 
\le \frac{4}{|y|} 
, \qquad 0<2|x|<|y| 
\label{eq: ineq A2}
\end{align}
for all $ \lambda>0 $. 
Using Fubini's theorem and integrating by parts, we have 
\begin{align}
h^{(2)}_q(v) 
=& \int_{(\lambda_0,\infty )} 
\kbra{ \int_{\lambda_0}^{\eta} (1-\cos \lambda v) \d \lambda } 
\Theta_q(\d \eta) 
\n \\
=& 
\int_{(\lambda_0,\infty )} 
\{ \eta-\lambda_0 - \varphi(\eta v) \eta + \varphi(\lambda_0 v) \lambda_0 \} 
\Theta_q(\d \eta) 
\n
\end{align}
for all $ v \in \bR $. 
Thus we see that 
\begin{align}
\frac{h^{(2)}_q(y-x) - h^{(2)}_q(y)}{x} 
=& \int_{(\lambda_0,\infty )} A^{(2)}_{\eta,y}(x) 
\eta \Theta_q(\d \eta) 
\n
\end{align}
where 
\begin{align}
A^{(2)}_{\eta,y}(x) 
= - \frac{\varphi(\eta (y-x)) - \varphi(\eta y)}{x} 
+ \frac{\lambda_0}{\eta} \cdot \frac{\varphi(\lambda_0 (y-x)) - \varphi(\lambda_0 y)}{x} . 
\n
\end{align}
By \eqref{eq: ineq A2}, we have 
\begin{align}
|A^{(2)}_{\eta,y}(x)| \le \frac{8}{|y|} 
, \qquad 0<2|x|<|y| 
\n
\end{align}
and hence, we obtain 
\begin{align}
\abra{ \frac{h^{(2)}_q(y-x) - h^{(2)}_q(y)}{x} } 
\le \frac{8}{|y|} 
\int_{(\lambda_0,\infty )} \eta \Theta_q(\d \eta) 
, \qquad 0<2|x|<|y| . 
\n
\end{align}
This implies \eqref{eq: ineq Cq y-x}. 

(ii) 
Note that 
\begin{align}
\frac{h^{(1)}_q(y+x) + h^{(1)}_q(y-x) - 2 h^{(1)}_q(y)}{x} 
= \int_0^{\lambda_0} \kbra{ A^{(1)}_{\lambda,y}(x) - A^{(1)}_{\lambda,y}(-x) } 
\frac{\lambda^2}{q+\theta(\lambda)} \d \lambda 
\n
\end{align}
and that 
\begin{align}
\frac{h^{(2)}_q(y+x) + h^{(2)}_q(y-x) - 2 h^{(2)}_q(y)}{x} 
= \int_{(\lambda_0,\infty )} \kbra{ A^{(2)}_{\eta,y}(x) - A^{(2)}_{\eta,y}(-x) } 
\eta \Theta_q(\d \eta) . 
\n
\end{align}
Since we have 
\begin{align}
\lim_{x \to 0} 
\kbra{ A^{(i)}_{\lambda,y}(x) - A^{(i)}_{\lambda,y}(-x) } 
= 0 
, \qquad i=1,2, 
\n
\end{align}
we obtain \eqref{eq: limit = 0} by the dominated convergence theorem. 

(iii) 
Since the estimate \eqref{eq: ineq A(1)} is valid for all $ x \neq 0 $ and $ y \in \bR $, 
we may apply the dominated convergence theorem to obtain 
\begin{align}
\lim_{x \to 0} \frac{1}{x} 
\int_{-\eps}^{\eps} \kbra{h^{(1)}_q(y+x) + h^{(1)}_q(y-x) - 2 h^{(1)}_q(y)} \d y = 0 . 
\label{eq: limit int = 0 (1)}
\end{align}

Note that 
\begin{align}
& \int_{-\eps}^{\eps} 
\kbra{ \varphi(\lambda (y+x)) + \varphi(\lambda (y-x)) - 2 \varphi(\lambda y) } \d y 
\n \\
=& 
\int_{-\eps}^{\eps} \d y 
\int_0^x (-\lambda) \{ \varphi'(\lambda (y-v)) - \varphi'(\lambda (y+v)) \} \d v 
\n \\
=& 
2 \int_0^x \{ \varphi(\lambda (\eps-v)) - \varphi(\lambda (\eps+v)) \} \d v 
\n
\end{align}
for all $ \lambda>0 $. 
Thus we have 
\begin{align}
\begin{split}
\int_{-\eps}^{\eps} \kbra{ A^{(2)}_{\eta,y}(x) - A^{(2)}_{\eta,y}(-x) } \d y 
= & - \frac{2}{x} 
\int_0^x \{ \varphi(\eta (\eps-v)) - \varphi(\eta (\eps+v)) \} \d v 
\\
&+ \frac{2 \lambda_0}{\eta x} 
\int_0^x \{ \varphi(\lambda_0 (\eps-v)) - \varphi(\lambda_0 (\eps+v)) \} \d v . 
\end{split}
\label{eq: int A(2)}
\end{align}
Since $ \varphi(v) $ is bounded in $ v \in \bR $, 
we see that the integral of the left hand side of \eqref{eq: int A(2)} 
is bounded in $ \eta>\lambda_0 $ and $ x \neq 0 $ 
and it converges to $ 0 $ as $ x \to 0 $. 
Since we have 
\begin{align}
\begin{split}
& \frac{1}{x} 
\int_{-\eps}^{\eps} \kbra{h^{(2)}_q(y+x) + h^{(2)}_q(y-x) - 2 h^{(2)}_q(y)} \d y 
\\
&= 
\int_{(\lambda_0,\infty )} \eta \Theta_q(\d \eta) 
\int_{-\eps}^{\eps} \kbra{ A^{(2)}_{\eta,y}(x) - A^{(2)}_{\eta,y}(-x) } \d y , 
\end{split}
\n
\end{align}
we see that this integral converges to $ 0 $ as $ x \to 0 $ 
by the dominated convergence theorem. 

Therefore we obtain \eqref{eq: limit int = 0}, which completes the proof. 
\qed
\end{proof}

\subsection{Limiting properties of the resolvent densities}

Define 
\begin{align}
p^{\sym}_t(x,y) 
= p^{\harm}_t(x,y) + p^{\harm}_t(-x,y) 
= p^{\harm}_t(x,y) + p^{\harm}_t(x,-y) 
\n
\end{align}
for $ t>0 $, $ x,y \in \bR $, and 
\begin{align}
u^{\sym}_q(x,y) 
= u^{\harm}_q(x,y) + u^{\harm}_q(-x,y) 
= u^{\harm}_q(x,y) + u^{\harm}_q(x,-y) 
\n
\end{align}
for $ q \ge 0 $, $ x,y \in \bR $. 
Then $ p^{\sym}_t(x,y) $ and $ u^{\sym}_q(x,y) $, respectively, are 
densities of the transition probability and the resolvent kernel, respectively, 
of the process $ \{ (|X_t|:t \ge 0),(P^{\harm}_x:x \in [0,\infty )) \} $ 
with respect to the symmetrizing measure $ h(y)^2 \d y $: 
\begin{align}
P^{\harm}_x(|X_t| \in A) = \int_A p^{\sym}_t(x,y) h(y)^2 \d y 
\n
\end{align}
for all $ t>0 $, $ x \in [0,\infty ) $ and $ A \in \cB((0,\infty )) $, 
and 
\begin{align}
P^{\harm}_x \ebra{ \int_0^{\infty } \e^{-qt} 1_A(|X_t|) \d t } 
= \int_A u^{\sym}_q(x,y) h(y)^2 \d y 
\n
\end{align}
for all $ q \ge 0 $, $ x \in [0,\infty ) $ and $ A \in \cB((0,\infty )) $.

\begin{lemma} \label{lem: h-trans}
Suppose that the conditions {\bf (A)} and {\bf (T)} are satisfied. Then 
the following assertions hold: 
\\ \quad {\rm (i)} 
For any $ q \ge 0 $ and $ y \in \bR \setminus \{ 0 \} $, 
\begin{align}
\lim_{x \to 0} u^{\sym}_q(x,y) = 
u^{\sym}_q(0,y) = 2 u^{\harm}_q(0,y) ; 
\label{eq: conti of hat u harm}
\end{align}
\quad {\rm (ii)} 
For any $ q>0 $, it holds that 
\begin{align}
\lim_{\eps \to 0+} \limsup_{x \to 0} \int_{-\eps}^{\eps} u^{\harm}_q(x,y) h(y) \d y = 0 . 
\label{pup limit2}
\end{align}
\end{lemma}

\begin{proof}%[Proof of Lemma \ref{lem: h-trans}]
(i) Let $ x,y \in \bR \setminus \{ 0 \} $. 
If $ q>0 $, we combine \eqref{eq: u harm qxy} with \eqref{eq: u harm qx} 
to cancel $ u_q(0) $, and then we have 
\begin{align}
u^{\harm}_q(x,y) 
= u^{\harm}_q(0,y) \cdot \frac{h_q(x)}{h(x)} 
- \frac{x}{h(x) h(y)} \frac{ h_q(y-x) - h_q(y) }{x} . 
\label{eq: u dagger q xy }
\end{align}
Then the identity \eqref{eq: u dagger q xy } is still valid for all $ q \ge 0 $. 
Hence we have 
\begin{align}
\begin{split}
u^{\sym}_q(x,y) 
=& 2 u^{\harm}_q(0,y) \cdot \frac{h_q(x)}{h(x)} 
\\
-& \frac{x}{h(x) h(y)} \cdot \frac{ h_q(y+x) + h_q(y-x) - 2 h_q(y) }{x} . 
\end{split}
\label{u never qh}
\end{align}
By Lemmas \ref{lem: bdd}, \ref{lem: conv to 0} and \ref{lem: estimate}, we obtain 
\eqref{eq: conti of hat u harm}. 

(ii) 
Integrating both sides of \eqref{u never qh} with respect to $ h(y) \d y $, we have 
\begin{align}
0 
\le& 
2 \int_{-\eps}^{\eps} u^{\harm}_q(x,y) h(y) \d y 
= 
\int_{-\eps}^{\eps} u^{\sym}_q(x,y) h(y) \d y 
\n \\
\le& \frac{2 h_q(x)}{h(x)} \int_{-\eps}^{\eps} u^{\harm}_q(0,y) h(y) \d y 
+ \frac{x}{h(x)} \int_{-\eps}^{\eps} \d y \frac{ u_q(y+x) + u_q(y-x) - 2 u_q(y) }{x} . 
\n
\end{align}
By Lemmas \ref{lem: bdd}, \ref{lem: conv to 0} and \ref{lem: estimate}, we obtain 
\begin{align}
\limsup_{x \to 0} \int_{-\eps}^{\eps} u^{\harm}_q(x,y) h(y) \d y 
\le \int_{-\eps}^{\eps} u^{\harm}_q(0,y) h(y) \d y . 
\n
\end{align}
Since 
\begin{align}
\int u^{\harm}_q(0,y) h(y) \d y 
= \frac{1}{u_q(0)} \int u_q(y) \d y 
= \frac{1}{qu_q(0)} < \infty , 
\n
\end{align}
we obtain \eqref{pup limit2}. 
\qed
\end{proof}

\begin{lemma} \label{lem: u dagger is conti at 0}
Suppose that the conditions {\bf (A)}, {\bf (B)} and {\bf (T)} are satisfied. 
Let $ q \ge 0 $ be fixed. 
Then it holds that 
\begin{align}
\lim_{x \to 0} u^{\harm}_q(x,y) = u^{\harm}_q(0,y) 
, \qquad y \in \bR \setminus \{ 0 \} . 
\label{eq: conti u harm}
\end{align}
Consequently, it holds that 
\begin{align}
\lim_{z \to 0} \frac{u_q^0(z,x)}{u_q^0(z,y)} = \frac{u^{\harm}_q(0,x)h(x)}{u^{\harm}_q(0,y)h(y)} 
, \qquad x,y \in \bR \setminus \{ 0 \} . 
\label{eq: conti u 0}
\end{align}
\end{lemma}

\begin{proof}
Let $ q \ge 0 $ and $ y \in \bR \setminus \{ 0 \} $ be fixed. 
Recall the identity \eqref{eq: u dagger q xy }: 
\begin{align}
u^{\harm}_q(x,y) 
= u^{\harm}_q(0,y) \cdot \frac{h_q(x)}{h(x)} 
- \frac{x}{h(x) h(y)} \frac{ h_q(y-x) - h_q(y) }{x} . 
\n
\end{align}
By Lemma \ref{lem: conv to 0}, (i) of Lemma \ref{lem: estimate} 
and the assumption {\bf (B)}, we obtain \eqref{eq: conti u harm}. 

Let $ q \ge 0 $ and $ x,y \in \bR \setminus \{ 0 \} $ be fixed. 
Then we obtain 
\begin{align}
\lim_{z \to 0} \frac{u_q^0(z,x)}{u_q^0(z,y)} 
= \lim_{z \to 0} \frac{u_q^{\harm}(z,x)h(x)}{u_q^{\harm}(z,y)h(y)} 
= \frac{u_q^{\harm}(0,x)h(x)}{u_q^{\harm}(0,y)h(y)} , 
\n
\end{align}
which proves \eqref{eq: conti u 0}. 
\qed
\end{proof}

\subsection{Transience of the $ h $-path process} \label{sec: transience}

Let us prove Theorem \ref{thm: trans of h-path}. 

%\begin{proof}[Proof of Theorem \ref{thm: trans of h-path}]
\noindent 
{\it Proof of Theorem \ref{thm: trans of h-path}} 
By a well-known theorem 
(see, e.g., \cite[Theorem 3.7.2]{MR2152573}), 
it suffices to prove the following: 
\\ \quad (i) 
The function 
\begin{align}
[0,\infty ) \ni x \mapsto \int_K u^{\sym}_0(x,y) h(y)^2 \d y 
\n
\end{align}
is lower-semicontinuous for any compact set $ K $ of $ [0,\infty ) $; 
\\ \quad (ii) 
There exists a nearly Borel function $ f $ which is positive almost everywhere 
such that 
\begin{align}
0 < \int_0^{\infty } f(y) u^{\sym}_0(x,y) h(y)^2 \d y < \infty . 
\label{eq: pos and fin}
\end{align}

Recall that 
\begin{align}
u^{\sym}_0(x,y) 
= \frac{2}{h(x)h(y)} \kbra{ h(x)+h(y) - \frac{h(x-y)+h(x+y)}{2} } 
- 2 \kappa . 
\n
\end{align}
The claim (i) is obvious by (i) of Lemma \ref{lem: h-trans} and by Fatou's lemma. 
The claim (ii) is also obvious; in fact, we may take 
$ f(y) = \min \{ 1,y^{-2} h(y)^{-2} \} $, which is a continuous function. 
Now the proof is complete. 
\qed
%\end{proof}

\subsection{The excursion measure of hitting a single point} \label{sec: hitting under vn}

Before closing this section, we give the following formula 
about the excursion measure of hitting a single point. 

\begin{theorem}
Suppose that the conditions {\bf (A)}, {\bf (B)} and {\bf (T)} are satisfied. 
Let $ a \in \bR \setminus \{ 0 \} $. 
Then it holds that 
\begin{align}
\vn ( T_{\{ a \}} < \zeta ) = \frac{1-\kappa h(a)}{h_B(a)} 
\n
\end{align}
where $ h_B(a) = 2h(a) - \kappa h(a)^2 $. 
\end{theorem}

\begin{proof}
Let $ x \in \bR $ and $ b \in \bR \setminus \{ a \} $. 
In our settings of symmetric L\'evy processes, 
Getoor's formula \cite[Theorem 6.5]{MR0185663} leads to 
\begin{align}
P_x(T_{\{ a \}} < T_{\{ b \}}) 
= \frac{h(a-b)-h(a-x)+h(b-x) - \kappa h(b-x) h(a-b)}{h_B(a-b)} . 
\n
\end{align}
Letting $ b=0 $ and using the symmetry $ h(-x) \equiv h(x) $, we have 
\begin{align}
P_x(T_{\{ a \}} < T_{\{ 0 \}}) 
= \frac{h(a)-h(a-x)+h(x) - \kappa h(x) h(a)}{h_B(a)} . 
\n
\end{align}
Let $ \eps>0 $. 
By the Markov property, we have 
\begin{align}
& \vn( \eps < T_{\{ a \}} < \zeta ) 
\label{eq: Hitting vn eq}
\\
=& \vn \ebra{ P^0_{X_{\eps}}( T_{\{ a \}} < \zeta ) ; \eps < T_{\{ a \}} \wedge \zeta } 
\n \\
=& \vn \ebra{ P_{X_{\eps}}( T_{\{ a \}} < T_{\{ 0 \}} ) ; \eps < T_{\{ a \}} \wedge \zeta } 
\n \\
=& \frac{1}{h_B(a)} \vn \ebra{ h(a)-h(a-X_{\eps})+h(X_{\eps}) - \kappa h(X_{\eps}) h(a) 
; \eps < T_{\{ a \}} \wedge \zeta } 
\n \\
=& \frac{1}{h_B(a)} P^{\harm}_0 \ebra{ 
\frac{X_{\eps}}{h(X_{\eps})} \cdot \frac{h(a)-h(a-X_{\eps})}{X_{\eps}} + 1 - \kappa h(a) 
; \eps < T_{\{ a \}} } . 
\label{eq: hitting of vn}
\end{align}
Now we let $ \eps \to 0+ $. 
On the one hand, using the assumption {\bf (B)}, Lemma \ref{lem: bdd} 
and (i) of Lemma \ref{lem: estimate}, 
we apply the dominated convergence theorem to see that 
the quantity \eqref{eq: hitting of vn} converges to 
\begin{align}
\frac{1 - \kappa h(a)}{h_B(a)} . 
\n
\end{align}
On the other hand, 
using monotone convergence theorem, we see that 
the quantity \eqref{eq: Hitting vn eq} converges to $ \vn(T_{\{ a \}} < \zeta) $, 
we obtain the desired result. 
\qed
\end{proof}

\section{Feller property of the $ h $-path process} \label{sec: Feller}

Define 
\begin{align}
T^{\harm}_tf(x) = P^{\harm}_x[f(X_t)] 
, \qquad t \ge 0 , f \in \cB_{+,b}(\bR) . 
\n
\end{align}
Then the Markov property implies that 
the family $ \{ T^{\harm}_t:t \ge 0 \} $ forms a transition semigroup: 
\\ \quad 
{\rm (T1)} 
$ T^{\harm}_{t+s} = T^{\harm}_t T^{\harm}_s $ for all $ t,s \ge 0 $; 
\\ \quad 
{\rm (T2)} 
$ T^{\harm}_0 $ equals the identity; 
\\ \quad 
{\rm (T3)} 
$ 0 \le f \le 1 $ implies that $ 0 \le T^{\harm}_tf \le 1 $ for all $ t \ge 0 $. 
\\
Note that {\rm (T3)} implies the contraction property: 
\\ \quad 
{\rm (T4)} 
$ \| T^{\harm}_tf \| \le \| f \| $ for $ t \ge 0 $ and $ f \in \cB_{+,b}(\bR) $. 
\\
We write the corresponding resolvent operator as 
\begin{align}
U^{\harm}_qf(x) = \int_0^{\infty } \e^{-qt} T^{\harm}_tf(x) \d t 
, \qquad q>0 , \ f \in \cB_{+,b}(\bR) . 
\label{eq: def of U dagger}
\end{align}
Then it is immediate that 
the family $ \{ U^{\harm}_q:q>0 \} $ satisfies the following properties: 
\\ \quad 
{\rm (U1)} 
$ U^{\harm}_q - U^{\harm}_r + (q-r) U^{\harm}_q U^{\harm}_r = 0 $ 
for $ q,r > 0 $; 
\\ \quad 
{\rm (U2)} 
$ 0 \le f \le 1 $ implies that $ 0 \le q U^{\harm}_q f \le 1 $ for $ q>0 $. 
\\
Note that {\rm (U2)} implies the {\em contraction property}: 
\\ \quad 
{\rm (U3)} 
$ \| q U^{\harm}_q f \| \le \| f \| $ for $ q>0 $ and $ f \in \cB_{+,b}(\bR) $. 
\\

\begin{lemma} \label{lem: excludes}
Suppose that the condition {\bf (A)} is satisfied. Then 
the condition {\bf (B)}, i.e., $ \lim_{x \to 0} \frac{x}{h(x)} = 0 $, 
implies $ v=0 $ in \eqref{def theta}. 
\end{lemma}

\begin{proof}
Suppose that $ v>0 $. Then we have 
$ \theta(\lambda) \ge v \lambda^2 $. 
Hence we obtain 
\begin{align}
h(x) 
\le \frac{1}{\pi} \int_0^{\infty } \frac{1-\cos \lambda x}{v \lambda^2} \d \lambda 
= \frac{|x|}{v \pi} C(2) = \frac{|x|}{2 v \pi} . 
\n
\end{align}
This prevents the condition {\bf (B)}. 
\qed
\end{proof}

Recall that the Feller property of the semigroup $ \{ T^{\harm}_t:t \ge 0 \} $ 
is stated precisely as follows: 
\\ \quad 
{\rm (F1)} 
$ T^{\harm}_t C_0(\bR) \subset C_0(\bR) $ for all $ t \ge 0 $; 
\\ \quad 
{\rm (F2)} 
$ \| T^{\harm}_t f - f \| \to 0 $ as $ t \to 0+ $ for all $ f \in C_0(\bR) $. 

In order to prove Theorem \ref{thm: Feller}, 
we shall prove the following 

\begin{proposition} \label{prop: Feller of Uq}
Suppose that the conditions {\bf (A)}, {\bf (B)} and {\bf (T)} are satisfied. Then 
the following statements hold: 
\\ \quad 
{\rm (i)} 
$ T^{\harm}_t f(x) \to f(x) $ as $ t \to 0+ $ for all $ x \in \bR $ and $ f \in C_0(\bR) $. 
\\ \quad 
{\rm (ii)} 
$ U^{\harm}_q C_0(\bR) \subset C_0(\bR) $ for each $ q>0 $. 
\end{proposition}

The proof of Proposition \ref{prop: Feller of Uq} will be given 
in Section \ref{sec: pf of Feller}. 
To deduce Theorem \ref{thm: Feller} from Proposition \ref{prop: Feller of Uq} 
is a kind of general argument, and so we omit it. 
See \cite[Proposition III.2.4]{MR1725357} for details.

\subsection{Feller property of the resolvent of the $ h $-path process} \label{sec: pf of Feller}

Now we are in a position to prove Proposition \ref{prop: Feller of Uq}. 

%\begin{proof}[Proof of Proposition \ref{prop: Feller of Uq}]
\noindent 
{\it Proof of Proposition \ref{prop: Feller of Uq}} 
(i) 
It is obvious since $ T^{\harm}_tf(x) = P^{\harm}_x[f(X_t)] $ 
and $ P^{\harm}_x $ is a probability measure on the c\`adl\`ag space $ \bD $. 

(ii) 
By the contraction property {\rm (U3)}, 
it suffices to show that $ U^{\harm}_q C_c(\bR) \subset C_0(\bR) $ 
where $ C_c(\bR) $ stands for the class of continuous functions $ \bR \to \bR $ with compact supports. 
Let $ f \in C_c(\bR) $ be fixed and let us prove that $ U^{\harm}_qf \in C_0(\bR) $. 

Recall that, for $ x \in \bR \setminus \{ 0 \} $, 
\begin{align}
U^{\harm}_qf(x) 
= \int f(y) u^{\harm}_q(x,y) h(y)^2 \d y 
= \frac{1}{h(x)} \int f(y) u^0_q(x,y) h(y) \d y . 
\label{eq: Uharmqf(x)}
\end{align}
Since $ f $ has compact support, and the functions $ h $ and $ f $ are continuous 
and $ u^{\harm}_q $ is continuous outside the origin, 
it is obvious that 
the function $ U^{\harm}_qf(x) $ is continuous at $ x \in \bR \setminus \{ 0 \} $. 

By \eqref{eq: Uharmqf(x)}, we have 
\begin{align}
|U^{\harm}_qf(x)| 
\le \frac{1}{h(x)} U^0_qf(x) \sup_{y \in \Supp(f)} |h(y)| . 
\n
\end{align}
As $ |x| \to \infty $, we have 
$ 1/h(x) \to \kappa < \infty $ ((iv) of Lemma \ref{lem: properties of hx}) 
and $ U^0_qf(x) \to 0 $. 
Since $ h $ is continuous ((i) of Lemma \ref{lem: properties of hx}), 
we see that $ U^{\harm}_qf(x) $ vanishes at infinity. 

Let us prove that 
the function $ U^{\harm}_qf(x) $ is continuous at $ x=0 $. 
Let $ \eps>0 $. Then, 
by Lemma \ref{lem: u dagger is conti at 0} 
and by the dominated convergence theorem, 
we obtain 
\begin{align}
\lim_{x \to 0} \int_{|x|>\eps} f(y) u^{\harm}_q(x,y) h(y)^2 \d y 
= \int_{|x|>\eps} f(y) u^{\harm}_q(0,y) h(y)^2 \d y . 
\n
\end{align}
We estimate the integral on the interval $ [-\eps,\eps] $ as 
\begin{align}
\abra{ \int_{-\eps}^{\eps} f(y) u^{\harm}_q(x,y) h(y)^2 \d y } 
\le \| fh \| \int_{-\eps}^{\eps} u^{\harm}_q(x,y) h(y) \d y . 
\n
\end{align}
Note that the right hand side coincide with 
\begin{align}
\| fh \| \int_{-\eps}^{\eps} u^{\sym}_q(x,y) h(y) \d y 
\n
\end{align}
by the symmetry $ h(-y) = h(y) $. 
Hence, by Lemma \ref{lem: h-trans}, we obtain 
\begin{align}
\lim_{\eps \to 0+} \limsup_{x \to 0} \abra{ \int_{-\eps}^{\eps} f(y) u^{\harm}_q(x,y) h(y)^2 \d y } 
= 0 . 
\n
\end{align}
Therefore we conclude that 
$ U^{\harm}_qf(x) \to U^{\harm}_qf(0) $ as $ x \to 0 $, 
which completes the proof. 
\qed
%\end{proof}

\subsection{Extremeness property}

Let us proceed to prove Corollary \ref{cor: extreme}. 

%\begin{proof}[Proof of Corollary \ref{cor: extreme}]
\noindent 
{\it Proof of Corollary \ref{cor: extreme}} 
By the Feller property of the semigroup $ \{ T^{\harm}_t:t \ge 0 \} $, 
we can prove, in the same way as Proposition \ref{prop: germ}, 
that the germ $ \sigma $-field $ \cF_{0+} $ 
is trivial under $ P^{\harm}_0 $. 
Since $ \vn $ is mutually absolutely continuous with respect to $ P^{\harm}_0 $, 
we see that 
the germ $ \sigma $-field $ \cF_{0+} $ 
is trivial also under $ \vn $. 
Hence, 
by Theorem \ref{thm: extremeness}, 
we conclude that $ \vn $ is an extreme direction. 
The proof is now complete. 
\qed
%\end{proof}

\subsection{Sample path behaviors}

Let us prove Corollaries 
\ref{cor: change signs of h-path} and 
\ref{cor: change signs of h-path-2}. 

%\begin{proof}
\noindent 
{\it Proof of Corollaries \ref{cor: change signs of h-path} 
and \ref{cor: change signs of h-path-2}} 
Set 
\begin{align}
\Omega^+_0 =& \kbra{ \text{$ \exists t_0>0 $ such that $ \forall t<t_0 $, $ X_t \ge 0 $} } , 
\n \\
\Omega^-_0 =& \kbra{ \text{$ \exists t_0>0 $ such that $ \forall t<t_0 $, $ X_t \le 0 $} } 
\intertext{and} 
\Omega^{+,-}_0 
=& \kbra{ \text{$ \exists \{ t_n \} $ with $ t_n \searrow 0 $ such that 
$ \forall n $, $ X_{t_n} X_{t_{n+1}} < 0 $} } . 
\n
\end{align}
Then the space $ \bD $ is decomposed into the disjoint union: 
\begin{align}
\bD = \Omega^+_0 \cup \Omega^-_0 \cup \Omega^{+,-}_0 . 
\n
\end{align}
Moreover, it is obvious that 
the three sets $ \Omega^+_0 $, $ \Omega^-_0 $ and $ \Omega^{+,-}_0 $ 
are all elements of $ \cF_{0+} $. 
Since $ \cF_{0+} $ is trivial under $ P^{\harm}_0 $, 
we see that 
only one of the three probabilities 
$ P^{\harm}_0(\Omega^+_0) $, $ P^{\harm}_0(\Omega^-_0) $ and $ P^{\harm}_0(\Omega^{+,-}_0) $ 
is one and the other two are zero. 
By the symmetry: $ P^{\harm}_0(X \in \cdot) = P^{\harm}_0(-X \in \cdot) $, 
we see that $ P^{\harm}_0(\Omega^+_0) $ and $ P^{\harm}_0(\Omega^-_0) $ coincide, 
which turn out to be zero. 
Therefore we conclude that $ P^{\harm}_0(\Omega^{+,-}_0)=1 $. 
This also proves that $ \vn((\Omega^{+,-}_0)^c)=0 $, 
which completes the proof. 
\qed
%\end{proof}

Now we prove Corollary \ref{cor: oscillation of stable}. 

%\begin{proof}
\noindent 
{\it Proof of Corollary \ref{cor: oscillation of stable}} 
Set 
\begin{align}
\Omega^+_{\infty } 
=& \kbra{ \text{$ \exists t_0>0 $ such that $ \forall t>t_0 $, $ X_t \ge 0 $} } , 
\n \\
\Omega^-_{\infty } 
=& \kbra{ \text{$ \exists t_0>0 $ such that $ \forall t>t_0 $, $ X_t \le 0 $} } 
\intertext{and} 
\Omega^{+,-}_{\infty } 
=& \kbra{ \text{$ \exists \{ t_n \} $ with $ t_n \nearrow \infty $ such that 
$ \forall n $, $ X_{t_n} X_{t_{n+1}} < 0 $} } . 
\n
\end{align}
Then the space $ \bD $ is decomposed into the disjoint union: 
\begin{align}
\bD = \Omega^+_{\infty } \cup \Omega^-_{\infty } \cup \Omega^{+,-}_{\infty } . 
\n
\end{align}
Since we have 
\begin{align}
P^{\harm}_0 \cbra{ \forall t>0 , \ X_t \neq 0 } = 1 
\n
\end{align}
by the local equivalence between $ P^{\harm}_0 $ and $ \vn $, 
we see that 
\begin{align}
\Omega^+_{\infty } \cup \Omega^-_{\infty } 
= 
\bigcup_{n=1}^{\infty } 
\kbra{ \text{$ \forall t>n $, $ X_n X_t > 0 $ } } 
\qquad \text{$ P^{\harm}_0 $-almost surely.} 
\n
\end{align}

Suppose that the process $ \{ (X_t),(P_x) \} $ 
is a symmetric stable process of index $ 1<\alpha <2 $. 
Then, from the original process $ \{ (X_t),(P_x) \} $, 
its $ h $-path process $ \{ (X_t),(P^{\harm}_0) \} $ 
inherits the scaling property: for any fixed $ c>0 $, 
\begin{align}
\cbra{ c^{-1/\alpha } X_{ct}:t \ge 0 } \law (X_t:t \ge 0) 
\qquad \text{under $ P^{\harm}_0 $.} 
\n
\end{align}
This implies that the probability 
\begin{align}
P^{\harm}_0 \cbra{ \text{$ \forall t>s $, $ X_s X_t > 0 $} } 
\n
\end{align}
for fixed $ s>0 $ does not depend on the choice of $ s>0 $. 
Hence we obtain 
\begin{align}
P^{\harm}_0 \cbra{ \Omega^+_{\infty } \cup \Omega^-_{\infty } } 
=& \lim_{n \to \infty } P^{\harm}_0 \cbra{ \text{$ \forall t>n $, $ X_n X_t > 0 $} } 
\n \\
=& \lim_{n \to \infty } P^{\harm}_0 \cbra{ \text{$ \forall t>1/n $, $ X_{1/n} X_t > 0 $} } 
\n \\
=& \lim_{n \to \infty } P^{\harm}_0 \cbra{ \text{$ X_t $ have the same sign for all $ t>0 $} } , 
\n
\end{align}
which proves to be zero by Corollary \ref{cor: change signs of h-path-2}. 
Hence we conclude that $ P^{\harm}_0(\Omega^{+,-}_{\infty })=1 $. 
By the transience of the $ h $-path process $ \{ (X_t),P^{\harm}_0 \} $ 
(Theorem \ref{thm: trans of h-path}), 
we have 
\begin{align}
\Omega^{+,-}_{\infty } = 
\kbra{ \limsup_{t \to \infty } X_t = \limsup_{t \to \infty } (-X_t) = \infty } 
\qquad \text{$ P^{\harm}_0 $-almost surely.} 
\n
\end{align}
Therefore the proof is complete. 
\qed
%\end{proof}

\subsection{Remark on a connection with a result of Ikeda--Watanabe} \label{sec: Ikeda-Watanabe}

Finally we make a remark on a connection with a result of Ikeda--Watanabe \cite{MR0451425}. 
Set 
\begin{align}
\Omega^{+,-}_1 
= \kbra{ \text{$ \exists \{ t_n \} $ with $ t_n \nearrow T_{\{ 0 \}} $ such that 
$ \forall n $, $ X_{t_n} X_{t_{n+1}} < 0 $} } . 
\n
\end{align}

\begin{theorem}[Theorem 3.3 of \cite{MR0451425}] \label{thm: Ikeda-Watanabe}
Suppose that, for any fixed $ q>0 $, 
\begin{align}
0 
< \liminf_{\eps \to 0+} \frac{u_q(0)-u_q(-\eps)}{u_q(0)-u_q(\eps)} 
\le \limsup_{\eps \to 0+} \frac{u_q(0)-u_q(-\eps)}{u_q(0)-u_q(\eps)} 
< \infty 
\label{eq: assump1}
\end{align}
and that 
\begin{align}
\lim_{\eps \to 0} \frac{u_q(x)-u_q(x + \eps)}{u_q(0)-u_q(\eps)} = 0 
, \qquad x \in \bR \setminus \{ 0 \} . 
\label{eq: assump2}
\end{align}
Then it holds that 
\begin{align}
P_x( \Omega^{+,-}_1 | T_{\{ 0 \}} < \infty ) = 1 
, \qquad x \in \bR \setminus \{ 0 \} . 
\label{eq: IW result}
\end{align}
\end{theorem}

Suppose that the conditions {\bf (A)}, {\bf (B)} and {\bf (T)} are satisfied. 
Then, by the symmetry $ u_q(x)=u_q(-x) $, 
we see that the assumption \eqref{eq: assump1} is satisfied. 
Since 
\begin{align}
\frac{u_q(x)-u_q(x + \eps)}{u_q(0)-u_q(\eps)} 
= \frac{\eps}{h(\eps)} \cdot \frac{h(\eps)}{h_q(\eps)} \cdot \frac{h_q(x+\eps)-h_q(x)}{\eps} , 
\n
\end{align}
we see, by Lemma \ref{lem: conv to 0}, (i) of Lemma \ref{lem: estimate} 
and the assumption {\bf (B)}, 
that the assumption \eqref{eq: assump2} is also satisfied. 
Hence we may apply Theorem \ref{thm: Ikeda-Watanabe} to obtain \eqref{eq: IW result}. 
The formula \eqref{eq: IW result} implies that 
\begin{align}
\vn((\Omega^{+,-}_1)^c \cap \{ \zeta < \infty \}) = 0 . 
\label{eq: Ikeda-Watanabe}
\end{align}
Through time reversal property (see \cite[Lemma 4.1]{MR2139028} and \cite[Lemma 5.2]{CFY}) 
of excursion paths with finite lifetime, 
the formula \eqref{eq: Ikeda-Watanabe} implies that 
\begin{align}
\vn \cbra{ (\Omega^{+,-}_0)^c \cap \{ \zeta < \infty \} } = 0 . 
\n
\end{align}
This is a special case of Corollary \ref{cor: change signs of h-path}.

{\bf Acknowledgements.} 
The author expresses his sincere gratitude to Professor Marc Yor 
for fruitful discussions about penalisation problems, 
from which the present study has originated. 
He expresses his hearty thanks to Professors Masatoshi Fukushima and Patrick J. Fitzsimmons, 
who were kind enough to draw his attension to their interesting papers 
\cite{CFY} and \cite{MR2247835}. 
He also thanks all of the participants of his seminar talks at Kyoto University 
for their valuable comments: 
Professors Shinzo Watanabe, Ichiro Shigekawa and Naomasa Ueki, 
and Doctors Nobuaki Sugimine and Yuko Yano.

\end{document}